\documentclass[12pt,oneside]{amsart}
\usepackage{mathtools,float}
\usepackage{cite}
\usepackage{color}
\usepackage[margin=1in]{geometry}
 \newtheorem{theorem}{Theorem}

 \newtheorem{lemma}[theorem]{Lemma}
 \newtheorem{corollary}{Corollary}
 
 \theoremstyle{remark}
 \newtheorem{remark}{Remark}
 \allowdisplaybreaks
\numberwithin{theorem}{section}
\numberwithin{corollary}{theorem}

\begin{document}

\title[Extremizers in the Rogosinski-Szeg\H{o} problem]
{Extremizers for the Rogosinski - Szeg\"o  estimate of the second coefficient in nonnegative sine polynomials}

\author[1]{Dmitriy Dmitrishin}
\address{Dmitriy Dmitrishin: Odessa National Polytechnic University, 1 Shevchenko Ave., Odesa 65044, Ukraine}
\email{dmitrishin@op.edu.ua}

\author[2]{Alexander Stokolos}
\address{Alexander Stokolos: Georgia Southern University, Statesboro GA, 30460, USA }
\email{astokolos@georgiasouthern.edu}

\author[3]{Walter Trebels}
\address{Walter Trebels: Department of Mathematics, AG Algebra, Technical University of Darmstadt, 64289 Darmstadt, Germany}
\email{trebels@mathematik.tu-darmstadt.de}

\maketitle

\begin{abstract}
For the class of sine polynomials $b_1\sin t+b_2\sin2t+...+b_N\sin Nt,\; (b_N\not= 0),$ which are nonnegative on $(0,\pi)$, W. Rogosinski and G. Szeg\"o derived, among other things, exact bounds for $|b_2|$ via the Luk\'acs
presentation of nonnegative algebraic polynomials and a variational type argument for exact bounds, but they did not find the extremizers. Within this algebraic framework, we construct explicit polynomials which attain these bounds and prove their uniqueness.
The proof uses the Fej\'er - Riesz representation of nonnegative trigonometric polynomials, a 7-band Toeplitz matrix of arbitrary finite dimension, and Chebyshev polynomials of the second kind and their derivatives.

\bigskip\noindent
{\sc Keywords}: Typically real polynomials, Chebyshev polynomials, extremal polynomials, nonnegative trigonometric polynomials.

\end{abstract}

\maketitle

\section{Introduction}
The classical problems of geometric complex analysis are related to the determination of the extremal properties of the functions $F(z)$ univalent in the unit disk $\mathbb{D}=\{z\in\mathbb{C}:|z|<1\}$ and having different normalizations. { The most popular one is the schlicht normalization: $F(0)=F'(0)-1=0$ (class $\mathcal{S}$)}. The significant drawback of the set of univalent functions is the lack of linearity: the sum of univalent functions is not necessarily a univalent function.   

 To avoid this drawback W.~Rogosinski in~\cite{WR32} introduced a class $\mathcal T$ of typically real functions. A holomorphic function  $F:\mathbb D\to\mathbb C$ belongs to $\mathcal T, F\in \mathcal T,$ if it satisfies 

i) $\forall z\in \mathbb R\cap\mathbb D,\; F(z)\in \mathbb R,$ 

ii) $\forall z\in \mathbb{D} \backslash\mathbb R,\; \text{Im}\{F(z)\}\cdot\text{Im}\{z\}>0.$ 

The class $\mathcal T$ has the convexity property, which appears to be convenient when solving various extremal problems. At the same time, many extremal estimates remain the same or are close to the analogous estimates in the class of univalent functions.

Denote by $\mathcal T_N\subset \mathcal T$ the set of typically real  polynomials   
$$
P(z)\in \mathcal T_N,\qquad P(z)=z+\sum_{j=2}^N a_j z^j,\qquad a_j\in \mathbb R.
$$
{Note, that $P(z)\in \mathcal T_N$ if and only if Im$(P(e^{it}))$ is a  sine polynomial nonnegative on $(0,\pi)$. 
\smallskip

In the pioneer work~\cite{RS50}, W.~Rogosinski and G.~Szeg\H{o} considered and discussed possible ways to solve a large variety of extremal problems for such polynomials. 
A very particular case of their results reads as follows: for $P(z)\in \mathcal T_N$ there holds the exact estimate}
\begin{equation}\label{eq:a2-estimate}
|a_2|\le\left\{
\begin{array}{ll}
2\mu_N,\;N\text{ is odd}, \\
2\eta_N,\;N\text{ is even},
\end{array}
\right.
\end{equation}
where $\mu_N=\cos\frac{2\pi}{N+3}$ is the largest root of the equation   $U_{(N+1)/2}(x)=0$, while {$\eta_N$ is the maximal root of $\displaystyle U^\prime_{\frac N2+1}(x)- U^\prime_{\frac N2}(x)=0.$} 
Here,  $U_j$ with $j\in\mathbb N_0$ denote the Chebyshev polynomials of the second kind and $U'_j$ their derivatives, {defined by
\begin{equation}\label{uj}    
U_j(x)=\frac{\sin(j+1)t}{\sin{t}}=\frac{z^{j+1}-z^{-j-1}}{z-z^{-1}},
\end{equation}
where $x=\cos t,$ and $z=e^{it}.$
} {The formulation in \eqref{eq:a2-estimate} is different than was originally written and is motivated by our deductions below. That both estimates coincide is shown in Remark \ref{3.1}.}

Let us note that writing the estimates in terms of the roots of Chebyshev polynomials or their derivatives is not only a technical matter but also gives a new conceptual insight. E.g., the occurrence of the derivative suggested feasible conjectures that led to a breakthrough in \cite{DSS22}. 

W.~Rogosinski and G.~Szeg\H{o} proved \eqref{eq:a2-estimate} by representing the trigonometric polynomial $\text{Im}\{P(e^{i\vartheta})\}$ through the Chebyshev polynomials of the second kind, used orthogonality properties with weights for these polynomials, and used the method of moments. By this technique, they obtained effectively all exact bounds for the coefficients $a_2$, $a_3$ as well as for $a_{N-1}$, $a_N.$ However, explicit extremizers and their uniqueness were not shown.{ We have found them below for the case $a_2$. 

The publication \cite{RS50} gave rise to several results by W.C.~Royster and T.~Suffridge in \cite{RS70} and in \cite{S83}. S.~Ruscheweyh \cite{R86}, using a remarkable theorem by O.~Sz\'asz \cite{Sz18}, rediscovered the estimates for $a_2$ and $a_3$, and expressed the bounds as the generalized eigenvalues for certain matrices. A table of numerical values of the coefficients (up to the degree ten) is published.  

Notable progress was made by D.K.~Dimitrov, C.A.~Merlo, and R.~Adreani in~\cite{DM02,AD03} using L. Fej\'er's method. In particular, extremizers for the quantities $a_{N-1}$, $a_N$, and $\text{Im}\{P(e^{i\vartheta})\}$ were constructed.} Furthermore, the exact upper and lower values of the quantity $a_{N-2}$ were found, and an extremizer was constructed for the case of odd $N$. It has been shown that in some problems, the extremizers are 
not uniquely determined. 

Classical problems of geometric complex analysis related to extremal stretching and contraction of the unit disk 
$\mathbb{D}$ by typically real polynomials were solved in~\cite{DSS22,DDS19,DSS23,MB89,CM71} (the extreme 
values and corresponding extremizers were found). Let us also note that various extremal problems in subclasses of 
typically real or univalent polynomials were considered, for example, in~\cite{RS70,S83,R86,IP21}.  

In \cite{RS50}, Rogosinski and Szeg\H{o} referred to an alternative way of solving extremal problems based on the Fej\'er-Riesz representation of a nonnegative trigonometric polynomial~\cite[6.5, Problem~41]{PS98} and a subsequent application of
the Rayleigh method for finding the extremum of the ratio of quadratic forms, which reduces to the problem of obtaining
the eigenvalues of a matrix pencil and their corresponding eigenvectors~\cite{MW68}. Yet, they \cite[p.115]{RS50} also remarked: ``In general, however, the method ... is not easily adaptable for obtaining explicit results, in particular when $N$ is large.'' {To prove our results, we take up this approach of converting the analytic problem into a linear algebraic one. However, this procedure is assuming big computational difficulties. For instance, its realization in \cite{DSS22} required computation of the eigenvalues and eigenvectors of arbitrarily size 5-band Toeplitz matrices. In the present article, we are working with 7-band matrices. 

Typically real polynomials $P(z)$ are of additional interest as a source of nonnegative trigonometric polynomials generated
by the nonnegative sine polynomials $\text{Im}\{P(e^{it})\}/\sin t$. Interesting applications of such polynomials in
approximation theory can be found, for instance, in~\cite{DM02,PT52}. 
Well-known are the Fej\'er, Gronwall-Jackson, and Egerv\'ary-Sz\'asz kernels~\cite{LF00,LF15,TG12,DJ11,ES28}. In a forthcoming paper, we will deal with this aspect.

{
\section{A brief survey over the main results and ideas of proof}
\subsection{Main results}\label{mainresults}

Our main results may briefly be outlined as follows.

\begin{itemize}    
   \item By Fej\'er's method we obtain \eqref{eq:a2-estimate} - see Corollaries \ref{cor:a2odd} and \ref{cor:a2even}.
    
    \item The extremal polynomials which attain $\max\{a_2\}$ and $\min\{a_2\}$ are unique - see Theorem \ref{ex-uni}.

    \item 
In the case of odd $N,$ the coefficients of the extremal polynomials $P^{odd}_{\max}(z)$ for $\max\{a_2\}$ are  given in Theorem \ref{tm:aj}, and a compact representation of $P^{odd}_{\max}(z)$ in Theorem \ref{tm:im-odd}. The resulting non-negative trigonometric polynomial
\begin{equation}\notag 
{\rm Im}\left(P^{odd}_{\max} (e^{it})\right)=\frac{1-\cos^2(\mu_N)}{N+3}\cdot\frac1{\sin t}\cdot\frac{\sin^2\frac{N+3}2t}{\left(\cos t - \mu_N\right)^2}
\end{equation}
is the unique extremizer for the estimate (0.7) in \cite{RS50} in the odd case.

\item 
In the case of even $N,$ the coefficients of the extremal polynomials $P^{even}_{\max}(z)$ for $\max\{a_2\}$ are  given in subsection \ref{sec:compform}, and a compact representation of  $P^{even}_{\max}(z)$ in Theorem \ref{tm:evencomp}. The resulting non-negative trigonometric polynomial 
{
\begin{equation}\notag 
    {\rm Im}\left(P^{even}_{\max}(e^{it})\right)=\frac{2(1-\eta_N^2)}{(N+2)(N+3)(N+4)}\cdot\frac1{1+\cos t}\cdot
\frac1{\sin t}\cdot\frac{\left(\frac{N+4}2\sin\frac{N+2}2t + \frac{N+2}2\sin\frac{N+4}2t\right)^2}{(\cos t-\eta_N)^2}.
\end{equation}} 
is a unique extremizer for the estimate (0.7) in \cite{RS50} in the even case.

\item     
The extremal polynomials for $\min\{a_2\}$ are
$$P^{odd}_{\min}(z)=-P^{odd}_{\max}(-z),\qquad P^{even}_{\min}(z)=-P^{even}_{\max}(-z).$$

\end{itemize}

\bigskip

{\bf Examples}. If $N=2,$ then 
$P^{even}_{\max}(z)=z+\frac12z^2.$ If $N=3,$ then $P^{odd}_{\max}(z)=z+z^2+\frac12z^3.$ Their derivation is given following Theorem \ref{tm:evencomp} and Theorem \ref{tm:im-odd}.

\subsection{Outline of the methods of proofs}

Because { the imaginary part of} a typically real polynomial on the unit circle  
is a non-negative sine polynomial on $[0,\pi]$ we are able to reduce the problem to a trigonometric one. Then, by factoring out $\sin t$, we can further reduce the problem to non-negative cosine polynomials with a simple relation between the original coefficients $a_j$ and
the cosine coefficients $\gamma_j$ given by formula \eqref{eq:a-to-gamma}. 

A core of the proof is the application of the Fej\'er-Riesz representation to the non-negative cosine polynomial which reduces the problem to the optimization of positive definite quadratic forms with coefficients $\delta_j.$ The relation between  $\gamma_j$ and $\delta_j$ - formula \eqref{eq:gamma-to-delta} - is more involved.  

The max/min problem for quadratic forms is reducible to finding the maximal/minimal eigenvalues of the corresponding matrix pencil. This leads to finding roots of the determinant of a specific 7-band Toeplitz matrix $\mathbf \Phi_N(x)$ which, in general, is a terribly complicated problem. Fortunately, we were able to solve it in our case, where we showed that the upper/lower bound in the Rogosinski-Szeg\"o estimate is the simple maximal/minimal eigenvalues of the matrix pencil. Thus, we regain the Rogosinski-Szeg\"o estimate \eqref{eq:a2-estimate} in \eqref{eq:modulus-a2-odd} and \eqref{eq:modulus-a2-even}.

Additionally, we get the existence and uniqueness of the extremizers. 

The next step is to determine the corresponding eigenvectors, {whose components $\delta_j$ will be denoted } $z^{(0)}_j(x)$ for the odd case and by $z^{(1)}_j(x)$ for the even case. This is done in Theorem \ref{tm:3} for odd $N$ and in Theorem \ref{tm:z1} for even $N.$ The situation is different for $N$ of different parity because we are looking for the maximal root of the Chebyshev polynomial in the case of odd $N$ and for the maximal root of its derivative in the case of even $N$.

Now we can specify $\gamma_j$ by \eqref{eq:gammas} and $a_j$ by \eqref {eq:as}.
In the odd case, it turns out, perhaps unsurprisingly, that the formulas for the coefficients $a_j$ of the extremal polynomial can be simplified, see Theorem \ref{tm:aj}. 

We point out that even though coefficient formulas are quite involved, especially in the even case, we were able to find the compact form for the extremal polynomials for both odd and even cases.
These extremal polynomials are represented as a sum of two rational functions, whose poles after summations turn into removable singularities. A similar approach has been used in \cite{MB89} for the solution of certain extremal problems.
The imaginary part of the extremal polynomials on the { unit circle} 
produces non-negative trigonometric polynomials \eqref{maxtrig} and \eqref{aseq:1}. 

\subsection{Notations}
Below, boldface letters will be used for matrices and vectors, i.e. $\mathbf \Phi_N(x)$ is a matrix, $\mathbf Z^{(0)}(x)$ is a vector, standard letters like $\Phi_N(x)$ and $ z_n^{(0)}(x)$ denote scalars. The subordered statements will be numbered by extending the numbering, e.g. a corollary of Theorem N will be denoted by Corollary N.1.  

\section{The estimate of $|a_2|$ via the Fej\'er-Riesz representation}\label{3.1}

\subsection{Transformation of the analytic problem into a linear algebraic one }

Let  $P(z)=z+\sum_{j=2}^N a_j z^j\in\mathcal T_N,$ hence Im$\{P(e^{it})\}\ge 0$ for $0\le t\le \pi.$ 
Following \cite{RS50} factor out the sine factor, i.e.  
\begin{equation}\label{eq:sinP}
    {\rm Im}\{P(e^{it})\}=(\sin{t})\left(\mathcal P(t)\right),\qquad \mathcal P(t)=\gamma_1+2\sum_{k=2}^N \gamma_k\cos(k-1)t. 
\end{equation}
 
$\mathcal P(t)$ is a non-negative cosine polynomial for $0\le t\le \pi,$ and
the coefficients $a_1,\ldots,a_N$ and $\gamma_1,\ldots,\gamma_N$ are related by the bijective relation
\begin{equation}\label{eq:a-to-gamma}
a_s=\gamma_s-\gamma_{s+2},\qquad s=1,\ldots,N.    
\end{equation} 
For convenience, in~\eqref{eq:a-to-gamma} we put $a_1=1$, $\gamma_{N+1}=\gamma_{N+2}=0.$ Therefore,
 $a_1=\gamma_1-\gamma_3=1$, $a_2=\gamma_2-\gamma_4$. 
By the Fej\'er-Riesz 
theorem the polynomial $\mathcal P(t)$  can be represented in the form 
\[
\mathcal P(t)=|\delta_1+\delta_2 e^{it}+\ldots+\delta_N e^{i(N-1)t}|^2,
\]
whence
\begin{equation}\label{eq:gamma-to-delta}
 \gamma_s=\sum_{j=1}^{N-s+1} \delta_j\delta_{j+s-1},\qquad s=1,\ldots,N.   
\end{equation}

Then 
\[
a_2=\gamma_2-\gamma_4=\sum_{j=1}^{N-1} \delta_j \delta_{j+1}-\sum_{j=1}^{N-3} \delta_j \delta_{j+3},
\;\;
1=\gamma_1-\gamma_3=\sum_{j=1}^N \delta_j^2-\sum_{j=1}^{N-2} \delta_j \delta_{j+2}.
\] 
Therefore,
\begin{align*}
\min\left\{\sum_{j=1}^{N-1}\delta_j\right.&\delta_{j+1}\left.-\sum_{j=1}^{N-3}\delta_j\delta_{j+3}:
\sum_{j=1}^N \delta_j^2-\sum_{j=1}^{N-2}\delta_j\delta_{j+2}=1\right\} \\
&\hspace{-2.26cm}\le a_2 \le
\max\left\{\sum_{j=1}^{N-1}\delta_j\delta_{j+1}-\sum_{j=1}^{N-3}\delta_j\delta_{j+3}:
\sum_{j=1}^N \delta_j^2-\sum_{j=1}^{N-2}\delta_j\delta_{j+2}=1\right\}.
\end{align*}

Associate symmetric matrices $\mathbf A$ and $\mathbf B$ of order $N\times N$ to the quadratic forms which are 7-band for $N\ge4$
\[
\sum_{j=1}^{N-1}\delta_j\delta_{j+1}-\sum_{j=1}^{N-3}\delta_j\delta_{j+3},
\quad
\mathbf A=
\begin{pmatrix}
0 & 1/2 & 0 & -1/2 & \ldots \\
1/2 & 0 & 1/2 & 0 & \ldots \\
0 & 1/2 & 0 & 1/2 & \ldots \\
-1/2 & 0 & 1/2 & 0 & \ldots \\
\ldots & \ldots & \ldots & \ldots & \ldots
\end{pmatrix};
\]
and
\[
\sum_{j=1}^N \delta_j^2-\sum_{j=1}^{N-2}\delta_j\delta_{j+2},\quad 
\mathbf  B=
\begin{pmatrix}
1 & 0 & -1/2 & 0 & \ldots \\
0 & 1 & 0 & -1/2 & \ldots \\
-1/2 & 0 & 1 & 0 & \ldots \\
0 & -1/2 & 0 & 1 & \ldots \\
\ldots & \ldots & \ldots & \ldots & \ldots
\end{pmatrix}.
\]
Note that the 7-band $N\times N$ matrix $\mathbf{A}$ has zeros along the main diagonal, the secondary diagonal above the main diagonal consists of $\frac{1}{2}$'s, the next secondary diagonal on top consists of zeros, the next one  on top of $-\frac{1}{2}$'s, the remaining secondary diagonals on top of zeros.  The same pattern occurs for the lower triangle. 

Similarly, $\mathbf{B}$ has 1's along the main diagonal, zeros on the first  secondary diagonal on top, then $-\frac{1}{2}$'s on the next  secondary diagonal on top, zeros on the remaining secondary diagonals on top. The same pattern occurs for the lower triangle.

Let $\lambda_1\le\ldots\le\lambda_N$ be the real roots of the equation $\det(\mathbf A-\lambda\mathbf  B)=0$ 
(note that $\mathbf A$ and $\mathbf B$ are symmetric matrices).
The numbers $\lambda_1,\ldots,\lambda_N$ are called the eigenvalues of the matrix pencil 
$\{\mathbf A-\lambda\mathbf B,\lambda\in\mathbb{C}\}$. 

The quadratic form $\sum_{j=1}^N \delta_j^2-\sum_{j=1}^{N-2} \delta_j \delta_{j+2}$ is positive definite
(see Lemma~\ref{le:A1}), therefore by the Rayleigh type theorem (see Theorem 4.2.2 in \cite[p.234]{HJ} and~\cite{MW68}) $\lambda_1\le a_2\le\lambda_N.$ To find the extremizers, it is necessary to know the eigenvectors corresponding to the eigenvalues $\lambda_1$ and
$\lambda_N$, that is, nontrivial solutions of the equations $(\mathbf A-\lambda_N\mathbf B)\mathbf Z=\mathbf 0$ and $(\mathbf A-\lambda_1\mathbf B)\mathbf Z=\mathbf 0$. 

Let the vector $\mathbf Z^{(0)}=(z_1^{(0)},\ldots,z_N^{(0)})^T$ be  an eigenvector of $\mathbf A-\lambda\mathbf B.$
Then, by formulas~\eqref{eq:gamma-to-delta} and \eqref{eq:a-to-gamma}, the coefficients of the
extremizer for the problem $\max\{a_2\}$ are defined for
$\ell=1,\ldots,N$ by
\begin{equation}\label{eq:gamma-0}
\gamma_\ell^{(0)}(\lambda)=\sum_{k=1}^{N-\ell+1} z_k^{(0)}(\lambda) z_{k+\ell-1}^{(0)}(\lambda),\qquad
a_\ell^{(0)}(\lambda)=\frac{\gamma_\ell^{(0)}(\lambda)-\gamma_{\ell+2}^{(0)}(\lambda)}{\gamma_1^{(0)}(\lambda)-\gamma_3^{(0)}(\lambda)},
\end{equation}
where we recall that $\gamma_{N+1}^{(0)}(\lambda)=\gamma_{N+2}^{(0)}(\lambda)=0$.
The coefficients of the extremizer for the problem $\min\{a_2\}$ are determined similarly.

Thus, the problem has been reduced to determining the eigenvalues of the matrix pencil $\{\mathbf A-\lambda\mathbf B,\lambda\in\mathbb{C}\}$ and their corresponding eigenvectors. It so happens that the cases of 
odd and even $N$ are fundamentally different, hence, they will be considered separately.

\subsection{Computation of the determinant $\bf \det(A-\lambda B)$}\label{sec:3.2}

Set $\lambda=2x$ and consider the matrix
\[
{\bf \Phi}_N (x)=2x\mathbf B-\mathbf A=
\begin{pmatrix}
2x & -1/2 & -x & 1/2 & 0 & 0 & \ldots \\
-1/2 & 2x & -1/2 & -x & 1/2 & 0 & \ldots \\
-x & -1/2 & 2x & -1/2 & -x & 1/2 & \ldots \\
1/2 & -x & -1/2 & 2x & -1/2 & -x & \ldots \\
0 & 1/2 & -x & -1/2 & 2x & -1/2 & \ldots \\
0 & 0 & 1/2 & -x & -1/2 & 2x & \ldots \\
\ldots & \ldots & \ldots & \ldots & \ldots & \ldots & \ldots  
\end{pmatrix}.
\]
We denote by $\Delta_N$ the determinant of this matrix. The determinants of the submatrices formed by discarding 
the first $k$ rows and $k$ columns in the original matrix will be denoted by $\Delta_{N-k}$, $k=1,\ldots,N-1$. By
Lemma~\ref{le:A2}, these determinants satisfy the {recursive} relation
\begin{equation}\label{eq:determ-relation}
\Delta_N-\sum_{j=1}^{10} 2^{-j} b_j\Delta_{N-2j}=0,\quad N\ge21,
\end{equation}  
{where the ${b_j}'$s are defined as in Lemma~\ref{le:A2}}. We will consider relation~\eqref{eq:determ-relation} as a linear
difference equation with constant coefficients of order 20. 
By $\{W_n\}_{n=1}^N$ we denote the solution of \eqref{eq:determ-relation} which satisfies
the initial conditions 
\begin{equation}\label{eq:initial-cond}
W_k=\Delta_k,\quad k=1,\ldots,20.
\end{equation} 

\subsubsection{\bf The case of odd $N$} 
It is shown in~\cite{DSS22,DSS23} that for the particular choices 
{$c_1=4x^2-3/2$, $c_2=-4x^4+2x^2-1/2$, $c_3=-1/2c_2$, $c_4=-1/8c_1$, $c_5=1/32$}  
the difference equation 
\begin{equation}\label{eq:diff-eq}
X_N-\sum_{j=1}^5 c_j X_{N-j}=0,\quad { N\ge 6,}
\end{equation}
has the solution 
$\Phi_N^{(1)}=2^{-N} U_{N+1}(x) U'_{N+1}(x)$
(for $ U_{N+1}$ see \eqref{uj}).

{The following observation is crucial. There holds 
\begin{equation}\label{eq:10-2}
\lambda^{10}-\sum_{j=1}^{10} b_j \lambda^{10-j}=\Bigg(\lambda^5-\sum_{j=1}^5 c_j\lambda^{5-j}\Bigg)^2.    
\end{equation}
This can be checked by direct computation and allows us to reduce the complicated equation \eqref{eq:determ-relation} to the known case of equation \eqref{eq:diff-eq}}.
Namely, if $\lambda$ is a single root of the characteristic equation $\lambda^5-\sum_{j=1}^5 c_j\lambda^{5-j}=0$ then this $\lambda$ is a root of multiplicity two of the left-hand side of 
\eqref{eq:10-2}. Thus, if $\lambda$ leads to the particular solution  $\Phi_N^{(1)}$ of \eqref{eq:diff-eq}, this implies the two particular solutions $\Phi_N^{(1)}$   and $\Phi_N^{(2)}=N\Phi_N^{(1)}$ of the equation $X_N-\sum_{j=1}^{10} b_j X_{N-j}=0$ (see e.g. \cite[p.76]{elaydi}). 

For odd $N$, the functions $\Phi_{\frac{N-1}2}^{(1)}$ and 
$\Phi_{\frac{N-1}2}^{(2)}$ are particular solutions of the equation $X_N-\sum_{j=1}^{10} b_j X_{N-2j}=0$. If we set
\[
Y_N^{(1)}=2^{-\frac{N}2}\Phi_{\frac{N-1}2}^{(1)}=2^\frac12 \frac1{2^N} U_{\frac{N+1}2}(x) U'_{\frac{N+1}2}(x),
\]
then $Y_N^{(1)}$ is a particular solution of \eqref{eq:determ-relation} since
\[
Y_N^{(1)}-\sum_{j=1}^{10} 2^{-j} b_j Y^{(1)}_{N-2j}=2^{-\frac{N}2}
\Bigg(\Phi_{\frac{N-1}2}^{(1)}-\sum_{j=1}^{10}2^{-j} b_j 2^j \Phi_{\frac{N-1}2-j}^{(1)}\Bigg)=0.
\] 
Similarly, the function
\[
Y_N^{(2)}=2^{-\frac{N}2}\Phi_{\frac{N-1}2}^{(2)}=
2^{-\frac12} \frac{N-1}{2^N} U_{\frac{N+1}2}(x) U'_{\frac{N+1}2}(x)
\]
is also a particular solution of equation~\eqref{eq:determ-relation}.
\begin{theorem}\label{tm:1} 
{ For odd  $N\in \mathbb N$ and $ -1< x< 1$
$$W_N:=\frac{N+3}{2^{N+2}} U_{\frac{N+1}2}(x) U'_{\frac{N+1}2}(x)
$$
is a solution of \eqref{eq:determ-relation} and \eqref{eq:initial-cond}. { The smallest (biggest) root of the equation $W_N(x)=0$ is simple. It is the minimal (maximal) eigenvalue of the matrix pencil $\{{\bf A-\lambda B}, \lambda\in \mathbb C\}$}.
}
\end{theorem}
\begin{proof}
The function $W_N$ is a linear combination of the
functions $Y_N^{(1)}$ and $Y_N^{(2)}$, hence this function satisfies \eqref{eq:determ-relation} when $N\ge21$. This
function also satisfies all relations~\eqref{eq:initial-cond}, which is verified by direct calculations. 
\end{proof}

By Theorem \ref{tm:1} we regain the Rogosinki-Szeg\"o estimate for odd $N$ in a natural way.

\begin{corollary}\label{cor:a2odd}
For typically real polynomials of odd degree $N$ there holds
\begin{equation}\label{eq:modulus-a2-odd}
|a_2|\le2\cos\frac{2\pi}{N+3}.
\end{equation}
\end{corollary}
\begin{proof}
The maximum and the minimum roots of the equation $W_N=0$ will be the maximum and the minimum roots of the 
equation $U_{\frac{N+1}2}(x)=0$, which are $\pm\cos\vartheta$, where $\sin\frac{N+3}2\vartheta=0$. Thus, 
$\vartheta=2\pi/(N+3)$, and $|a_2|\le2\cos\vartheta=2\cos\frac{2\pi}{N+3}$.  
\end{proof}

\subsubsection{\bf The case of even $N$}
\begin{theorem}\label{tm:2} 
{For even $N\in \mathbb N$ and $ -1< x< 1$
$$V_N=\frac1{2^{N+2}} \left( \left(U'_{\frac{N}2+1}(x)\right)^2-\left(U'_{\frac{N}2}(x)\right)^2\right)
$$
is a solution of \eqref{eq:determ-relation} and \eqref{eq:initial-cond}.
{The smallest (biggest) root of the equation $V_N(x)=0$ is simple. It is the minimal (maximal) eigenvalue of the matrix pencil $\{{\bf A-\lambda B}, \lambda\in \mathbb C\}$}.
}
\end{theorem}
\begin{proof}
Let $k=N/2$. By Lemma~\ref{le:A4},
\begin{equation}\label{eq:deltan}
  \Delta_N=\Delta_{2k}=\frac1{2^{2k+2}}\frac1{1-x^2}
\Big((k+2)^2(U_k(x))^2-(k+1)^2(U_{k+1}(x))^2\Big).  
\end{equation}
{Using \eqref{uj} with $z=e^{it}$ and $x=\cos t,$ we obtain} 
\[
\Delta_{2k}=\frac{{-}1}{2^{{ 2}k}}\frac{z^4}{(1-z^2)^4}\Big((k+2)^2(z^{k+1}-z^{-k-1})^2-(k+1)^2(z^{k+2}-z^{-k-2})^2\Big).
\] 
Set
$R_k=z^{-4}(1-z^2)^4\Delta_{2k}$ and rewrite the preceding equation in terms of $R_k$ to obtain 
\[
R_k=4^{-k}\Big(z^2 R_k^{(1)}+z^{-2} R_k^{(2)}-z^4 R_k^{(3)}-z^{-4}R_k^{(4)}-2R_k^{(5)}\Big), 
\]
where
$$R_k^{(1)}=(k+2)^2z^{2k},\quad R_k^{(2)}=(k+2)^2 z^{-2k},\quad R_k^{(3)}=(k+1)^2 z^{2k},$$ 
$$R_k^{(4)}=(k+1)^2 z^{-2k},\quad R_k^{(5)}=(2k+3).$$

Now, substitute $\Delta_k$ by $R_k$ in the left-hand side of \eqref{eq:determ-relation}, take into account that $x=\frac12(z+z^{-1})$ i.e. replace $b_j$ by $\hat b_j$ (Lemma \ref{le:A3}) to arrive at the linear expression (in $R_j$)
\begin{align*}
R_k-\sum_{j=1}^{10} 2^{-j} \hat{b}_j &R_{k-j}=4^{-k} 
\bigg[\Big(z^2 R_k^{(1)}+z^{-2}R_k^{(2)}-z^4 R_k^{(3)}-z^{-4}R_k^{(4)}-2R_k^{(5)}\Big) \\
&\hspace{-1cm}-\sum_{j=1}^{10} 2^{-j} \hat{b}_j
\Big(z^2 R_{k-j}^{(1)}+z^{-2}R_{k-j}^{(2)}-z^4 R_{k-j}^{(3)}-z^{-4}R_{k-j}^{(4)}-2R_{k-j}^{(5)}\Big)\bigg].
\end{align*}
Each function $R_k^{(s)}$ ($1\le s\le 5$) satisfies equation~\eqref{eq:psi}. This yields
$ R_k-\sum_{j=1}^{10}2^{-j}\hat{b}_j R_{k-j}=0$, 
 hence the same equation is valid for $\Delta_{2k},$ which by the third equation in  Lemma~\ref{le:A4}, statement b),
implies that the function 
\[
V_N=\frac1{2^{N+2}} \left( \left(U'_{\frac{N}2+1}(x)\right)^2-\left(U'_{\frac{N}2}(x)\right)^2\right)
\]
satisfies~\eqref{eq:determ-relation}. This function also satisfies all relations~\eqref{eq:initial-cond}, which can be checked by direct calculations. 
\end{proof}
Thus, by Theorem \ref{tm:2}, the Rogosinki-Szeg\"o estimate for even $N$ is regained.

\begin{corollary}\label{cor:a2even}
For typically real polynomials of even degree $N$ there holds
\begin{equation}\label{eq:modulus-a2-even}
|a_2|\le2(1-2\nu_N^2),
\end{equation}
where $\nu_N$ is the smallest positive root of the equation $U'_{N+2}(x)=0$.
\end{corollary}
\begin{proof}
By Lemma~\ref{le:A4} e),
\begin{equation}\label{eq:delNeven}
\Delta_N=\frac1{2^{N+5}}\frac{(-1)^\frac{N}2}{\sqrt{1-x^2}}
U'_{N+2}\Bigg(\sqrt{\frac{1+x}2}\Bigg)U'_{N+2}\Bigg(\sqrt{\frac{1-x}2}\Bigg).
\end{equation}

Let $\nu_N^{\max}$ and $\nu_N^{\min}$ be the largest and the smallest positive roots, respectively, of the equation $U'_{N+2}(x)=0$. Then, the largest root of the equation $\Delta_N=0$ does not exceed the value $\max\{2(\nu_n^{\max})^2-1,1-2(\nu_N^{\min})^2\}$. By Lemma \ref{le:A9},
$(\nu_N^{\min})^2+(\nu_N^{max})^2<1,$ whence $2(\nu_N^{\max})^2-1<1-2(\nu_N^{\min})^2$. Thus, this largest root equals $1-2(\nu_N^{\min})^2$, which 
proves the statement.

\end{proof}

{\sc Remark \ref{3.1}.} The estimates \eqref{eq:modulus-a2-odd} and \eqref{eq:modulus-a2-even} coincide with those originally given by Rogosinski and Szeg\"o
$$
|a_2|\le\left\{
\begin{array}{ll}
2\cos\frac{2\pi}{N+3},\;N\text{ is odd}, \\
2\cos \theta,\;\; \;N\text{ is even},
\end{array}
\right.
$$
where $\theta$ is the smallest positive root of the equation
\begin{equation}\label{eq:RSt}
(N+4)\sin\frac{N+2}2\vartheta+(N+2)\sin\frac{N+4}2\vartheta=0.
\end{equation}
For the case of odd $N$ the coincidence is obvious. Concerning the case of even $N$, observe that formula a) in Lemma \ref{le:A4} implies the following: the equation $U'_{N+2}(x)=0$ is equivalent to 
$$(N+4)U_{N+1}(x)-(N+2)U_{N+3}(x)=0.$$  
Let $\nu=\cos\tau$ be a root of this equation then, having in mind the definition of the Chebyshev polynomials \eqref{uj}, we arrive at
\begin{equation}\label{eq:tau}
    (N+4)\sin(N+2)\tau-(N+2)\sin(N+4)\tau=0.
\end{equation}
Now, the right hand side of \eqref{eq:modulus-a2-even}  can be written as $1-2\nu^2=-\cos(2\tau)=\cos(2\tau+\pi)=\cos \vartheta,$ therefore \eqref{eq:tau} in terms of $\theta$ yields \eqref{eq:RSt}. 

The quantities $\mu_N$ (largest root of $U_{\frac{N+1}2}(x)=0$) and $\eta_N=1-2\nu^2$ (largest root of $U'_{\frac{N}2+1}(x) - U'_{\frac{N}2}(x)=0$) have been introduced to emphasize the structural uniformity of the estimate \eqref{eq:a2-estimate}.
}

\begin{theorem}\label{ex-uni}
For typically real polynomials there exist unique extremal polynomials attaining the estimates in Corollaries \ref{cor:a2odd} and \ref{cor:a2even}.  
\end{theorem}

\begin{proof}
The existence follows from the existence of a solution for the quadratic form Rayleigh-type extremal problem - see \cite{R86} for the derivation in similar problems.
The uniqueness follows from the fact that the maximum and the minimum eigenvalue of the matrix pencil $\{A-\lambda B,\lambda\in\mathbb{C}\}$ divided by 2 are simple roots of the equation 
$\displaystyle U_{\frac{N+1}2}(x)U'_{\frac{N+1}2}(x)=0$ (see Theorem \ref{tm:3}), or the equation
\begin{equation}\label{eq:U'U'}
    \frac1{\sqrt{1-x^2}}U'_{N+2}\left(\sqrt{\frac{1+x}2}\right)U'_{N+2}\left(\sqrt{\frac{1-x}2}\right)=0,
\end{equation}
due to \eqref{eq:delNeven}.
\end{proof}

\section{Explicit formulas for the extremizers and their coefficients in the case of odd $N.$}
\subsection{Eigenvectors of the matrix pencil $\bf \{A-\lambda B,\lambda\in\mathbb{C}\}$}
\begin{theorem}\label{tm:3}
The solution of the system of linear equations
\[
\left(2\cos j\frac{2\pi}{N+3}\mathbf B-\mathbf A\right)\mathbf Z=\mathbf 0,\quad j=1,\ldots,{\frac{N+1}2},
\]
is the one-parameter family 
\[c\mathbf Z^{(0)}\left(\cos{j\frac{2\pi}{N+3}}\right),\quad  \mathbf Z^{(0)}(x)=(z_1^{(0)}(x),\ldots,z_N^{(0)}(x))^T,\]
where $c\in\mathbb{R},\; z_N^{(0)}=1,$ and for $k=1,\ldots,(N-1)/2$
$$
z_{2k-1}^{(0)}(x)=U_{k-1}(x)U_{k-1}(x),\quad z_{2k}^{(0)}(x)=U_{k-1}(x)U_k(x).$$ 
The symbol $T$ denotes transposition. 
\end{theorem}\label{tm:4}

\begin{proof}[Proof\nopunct] 
follows from Lemma~\ref{le:A5}.
\end{proof}
{The following useful property of the eigenvectors can be easily verified from the definition of $U_j(x)$}: 
\begin{equation}\label{eq:zz}
z_k^{(0)}\left(\cos{j\frac{2\pi}{N+3}}\right)=z_{N-k+1}^{(0)}\left(\cos{j\frac{2\pi}{N+3}}\right),\quad
k=1,\ldots,N.
\end{equation}
In what follows, we will direct our attention to the maximum eigenvalue $\mu_N=\cos\left(2\pi/(N+3)\right)$, and its corresponding eigenvector $\mathbf Z^{(0)}(\mu_N)$.\\

\subsection{\bf Computing the coefficients of extremizers}

\begin{theorem}\label{tm:aj} Let $N$ be an odd integer and $\mu_N=\cos\frac{2\pi}{N+3}.$ Then, the coefficients in formula \eqref{eq:gamma-0} for $\displaystyle j=1,\ldots,(N-1)/2$ are explicitly given by
\begin{align}\label{eq:a_2j+1}
a_{2j+1}&=\frac2{N+3}\Bigg[\mu_N U_{j-1}(\mu_N)U_j(\mu_N)+\left(\frac{N+3}2-j\right)(U_j(\mu_N))^2 \nonumber \\
&\hspace{4.5cm} -\left(\frac{N+1}2-j\right)(U_{j-1}(\mu_N))^2\Bigg],
\end{align}
\begin{align}\label{eq:a_2j}
a_{2j}&=\frac4{N+3}U_{j-1}(\mu_N)\Bigg[\left(\frac{N+3}2-j\right)U_j(\mu_N) \nonumber \\
&\hspace{5cm} -\mu_N\left(\frac{N+1}2-j\right)U_{j-1}(\mu_N)\Bigg].
\end{align}    
\end{theorem}

\begin{proof} To determine the coefficients $a_j,$ by \eqref{eq:gamma-0} we first have to calculate the $\gamma$'s for the particular case $\lambda=\mu_N.$ To keep the formulas compact, we will omit the argument $\mu_N$ in the computations. Thus, by Theorem \ref{tm:3}, for $j=1,\ldots,(N+1)/2$ we have  
\begin{align*}
\gamma_{2j-1}(\mu_N)&=\sum_{k=1}^{N-2j+2} z_k^{(0)} z_{k+2j-2}^{(0)} 
=\sum_{k=1}^{\frac{N-2j+3}2} \left( z_{2k}^{(0)} z_{2(k+j-1)}^{(0)}
+ z_{2k-1}^{(0)} z_{2k+2j-3}^{(0)}\right) \\
&\begin{aligned}
=\sum_{k=1}^{\frac{N-2j+3}2}
\bigl(U_{k-1}U_kU_{k+j-2}U_{k+j-1}+(U_{k-1})^2(U_{k+j-2})^2\bigr);& 
\end{aligned} \\
\end{align*}
\begin{align}\label{eq:gamma2j}
\gamma_{2j}(\mu_N)&=\sum_{k=1}^{N-2j+1} z_k^{(0)} z_{k+2j-1}^{(0)}=\sum_{k=1}^{\frac{N-2j+1}2}\left( z_{2k}^{(0)} z_{2k+2j-1}^{(0)}
+ z_{2k-1}^{(0)} z_{2k+2j-2}^{(0)}\right)\notag \\
&\begin{aligned}
=\sum_{k=1}^{\frac{N-2j+1}2}
U_{k-1}U_k(U_{k+j-1})^2+
\sum_{k=1}^{\frac{N-2j+1}2} (U_{k-1})^2 U_{k+j-2}U_{k+j-1}.& \\
\end{aligned}
\end{align} 

We start with computing 
$$\gamma_{2j}(\mu_N)-\gamma_{2j+2}(\mu_N)=: I_{j,1}+I_{j,2},$$
where we rearrange the contributing sums in the following way:

\begin{align*}
&I_{j,1}:=\sum_{k=1}^{\frac{N+1}2-j}(U_{k-1})^2U_{k+j-2}U_{k+j-1}
-\sum_{k=1}^{\frac{N-1}2-j}U_{k-1}U_k(U_{k+j})^2, \\
&I_{j,2}:=\sum_{k=1}^{\frac{N+1}2-j}U_{k-1}U_{k}(U_{k+j-1})^2
-\sum_{k=1}^{\frac{N-1}2-j}(U_{k-1})^2U_{k+j-1}U_{k+j}. \\
\end{align*}

The use of the relation 
\begin{equation}\label{eq:UU}  
U_{k-1}(x)U_{k+j}(x)=U_k(x)U_{k+j-1}(x)-U_{j-1}(x)
\end{equation}
leads to 
\begin{align*}
I_{j,1}=(U_0)^2U_jU_{j-1}&+\sum_{k=1}^{\frac{N-1}2-j}\left[(U_k)^2U_{k+j-1}U_{k+j}
-(U_k)^2U_{k+j-1}U_{k+j}\right]& \\
\\
&+U_{j-1}\sum_{k=1}^{\frac{N-1}2-j}U_kU_{k+j}=U_jU_{j-1}+U_{j-1}\sum_{k=1}^{\frac{N-1}2-j}U_kU_{k+j},\\
\end{align*}
(This equation even holds for all $x$). Again, using \eqref{eq:UU} we obtain

\begin{align*}
I_{j,2}&=U_{\frac{N-1}2-j}U_{\frac{N+1}2-j}\left(U_{\frac{N+1}2-1}\right)^2 + \sum_{k=1}^{\frac{N+1}2-j} \left[U_{k-1}U_{k}(U_{k+j-1})^2-U_{k-1}U_{k}(U_{k+j-1})^2\right.\\
&\left.+U_{j-1}U_{k-1}U_{k+j-1}\right]=U_{j-1}U_{j}+U_{j-1}\sum_{k=1}^{\frac{N-1}2-j}U_{k-1}U_{k+j-1}.
\end{align*}
Here, the last equality holds by \eqref{eq:zz} for $j=1,...,(N-3)/2.$ Summarizing, we arrive at 
\begin{align*}
\gamma_{2j}-\gamma_{2j+2}&=2U_{j-1}U_j
+U_{j-1}\left(\sum_{k=1}^{\frac{N-1}2-j}U_kU_{k+j}+\sum_{k=1}^{\frac{N-1}2-j}U_{k-1}U_{k+j-1}\right) \\
&=2U_{j-1}U_j+2U_{j-1}\sum_{k=1}^{\frac{N-1}2-j}U_kU_{k+j},\qquad (j=1,...,(N-3)/2).
\end{align*}.

We handle the case of the $\gamma$'s with odd index analogously. First, observe that, by Lemma \ref{le:A7} (for $n=(N-1)/2$),
\[
\gamma_1(\mu_N)-\gamma_3(\mu_N)=1+\sum_{k=1}^{\frac{N+1}2-1} (U_k(\mu_N))^2=\frac{N+3}{4\sin^2\frac{2\pi}{N+3}}.
\]
There remains to consider for $j=1,...,\frac{N+1}2$
$$
\gamma_{2j+1}(\mu_N)-\gamma_{2j+3}(\mu_N)=:II_{j,1}+ II_{j,2}
$$
where, on account of Theorem \ref{tm:3},

\begin{align*}
&II_{j,1}:= \sum_{k=1}^{\frac{N+1}2-j}(U_{k-1})^2(U_{k+j-1})^2
-\sum_{k=1}^{\frac{N-1}2-j}U_{k-1}U_kU_{k+j}U_{k+j+1}, \\
&II_{j,2}:= \sum_{k=1}^{\frac{N+1}2-j}U_{k-1}U_kU_{k+j-1}U_{k+j}
-\sum_{k=1}^{\frac{N-1}2-j}(U_{k-1})^2(U_{k+j})^2.
\end{align*}{ Using the original relation  \eqref{eq:UU}, and with $j$ replaced by $j+1$ in one summation, we get
}
\begin{align*}
II_{j,1}&=\sum_{k=1}^{\frac{N+1}2-j}(U_{k-1})^2(U_{k+j-1})^2
-\sum_{k=1}^{\frac{N-1}2-j}(U_{k})^2(U_{k+j})^2
+U_j\sum_{k=1}^{\frac{N-1}2-j}U_kU_{k+j} \\
&=(U_0)^2(U_j)^2+U_j\sum_{k=1}^{\frac{N-1}2-j}U_kU_{k+j},\\
II_{j,12}&=\sum_{k=1}^{\frac{N+1}2-j}U_{k-1}U_kU_{k+j-1}U_{k+j}
-\sum_{k=1}^{\frac{N-1}2-j}U_{k-1}U_kU_{k+j-1}U_{k+j}\\
&+U_{j-1}\sum_{k=1}^{\frac{N-1}2-j}U_{k-1}U_{k+j}=U_{j-1}\sum_{k=1}^{\frac{N-1}2-j}U_{k-1}U_{k+j}.\\
\end{align*}

Here, we take into account that 
\[
U_{\frac{N+1}2-j}(\mu_N)U_{\frac{N+1}2-1}(\mu_N)U_{\frac{N-1}2-j}(\mu_N)U_{\frac{N+1}2}(\mu_N)=0.
\]
Then, 
\[
\gamma_{2j+1}(\mu_N)-\gamma_{2j+3}(\mu_N)=(U_j)^2+U_j\sum_{k=1}^{\frac{N-1}2-j}U_{k}U_{k+j}
+U_{j-1}\sum_{k=1}^{\frac{N-1}2-j}U_{k-1}U_{k+j}.
\]

Thus, the following formulas are obtained for the extremizer coefficients, {where, until the end of the proof, the value of the Chebyshev polynomials are computed at $\mu_N,$ i.e. $U_j$ means $U_j(\mu_N)$} : 
\begin{align*}
a_{2j+1}&=\frac4{N+3}\left(\sin^2\frac{2\pi}{N+3}\right)\Bigg[(U_j)^2
+U_j\sum_{k=1}^{\frac{N-1}2-j}U_kU_{k+j} 
+U_{j-1}\sum_{k=1}^{\frac{N-1}2-j}U_{k-1}U_{k+j}\Bigg],\\
a_N&=\frac4{N+3}\sin^2\frac{2\pi}{N+3},\hspace{8.5cm} j=0,\ldots,\frac{N-3}2.\\
a_{2j}&=\frac4{N+3}\left(\sin^2\frac{2\pi}{N+3}\right)\Bigg[2U_{j-1}U_j
+2U_{j-1}\sum_{k=1}^{\frac{N-1}2-j}U_{k}U_{k+j}\Bigg],\hspace{2.6cm} j=1,\ldots,\frac{N-1}2.
\end{align*}

To simplify the preceding formulas, by Lemma~\ref{le:A8}, we obtain
\begin{align*}
a_{2j+1}=&\frac{2U_j}{N+3}\Bigg[\left(\frac{N-3}2-j\right)\mu_N U_{j-1}-\left(\frac{N-1}2-j\right)U_{j-2}+2\mu_N^2 U_j\Bigg]\\
&\hspace{2.6cm}+\frac{4(1-\mu_N^2)U_j^2}{N+3}+\frac{2U_{j-1}}{N+3}\Bigg[\left(\frac{N+3}2-j\right)\mu_N U_j-\left(\frac{N+1}2-j\right)U_{j-1}\Bigg]\\
&\hspace{-.25cm} ={\frac{2}{N+3}\Bigg[\left(\frac{N-3}2-j\right)\mu_N U_jU_{j-1}-\left(\frac{N-1}2-j\right) U_j(2\mu_N U_{j-1} -U_j)}\\
&{\hspace{1.8cm} +2\mu_N^2 U_j^2 + 2(1-\mu_N^2)U_j^2 +\mu_N U_{j-1}U_j\left(\frac{N+3}2-j\right) -U_{j-1}^2\left(\frac{N+1}2-j\right) \Bigg].}\\
\end{align*}
For the last equality, we use the recurrence formula for the Chebyshev polynomials, $U_j(x)=2xU_{j-1}(x)-U_{j-2}(x).$ 
Now, combine the coefficients belonging to $ U_j U_{j-1}$, then those belonging to $U_j^2,$ then to $U_{j-1}^2,$ which yields \eqref{eq:a_2j+1} for $ 1 \le j\le (N-1)/2 .$ Formula \eqref{eq:a_2j} is derived similarly.
\end{proof}

\subsection{\bf Compact form for the extremizers} 

The proof of the following theorem can be obtained from formulas \eqref{eq:a_2j+1} and \eqref{eq:a_2j} by summing up a geometric progression and its derivative. However, to save the reader's time, let us omit the computations and just prove that the final formula provides the desired extremizers. 

\begin{theorem}\label{tm:im-odd}  The following representations of $P^{odd}_{\max}(z)$ holds:
$$
P^{odd}_{\max}(z)=P^o_1(z)+P^o_2(z),
$$
where
\begin{align}\label{ext1}
 P^o_1(z)&=\frac z{1-2z\cos\frac{2\pi}{N+3}+z^2}, \\
P^o_2(z)&=\frac4{N+3}\cdot\sin^2\frac{2\pi}{N+3}\cdot \frac{z^3}{1-z^2}\cdot\frac{1-z^{N+3}}{\left(1-2z\cos\frac{2\pi}{N+3}+z^2\right)^2}.\label{ext}%
\end{align}
Hence, the resulting { non-negative trigonometric polynomial} has the form
\begin{equation}\label{maxtrig}
{\rm Im}\left(P^{odd}_{\max} (e^{it})\right)=\frac{\sin^2\frac{2\pi}{N+3}}{N+3}\cdot\frac1{\sin t}\cdot\frac{\sin^2\frac{N+3}2t}{\left(\cos t - \cos\frac{2\pi}{N+3}\right)^2}.
\end{equation}
For the problem $\min\{a_2\}$, the extremizer is obtained from the maximum extremizer by alternating signs for even powers, i.e. $P^{odd}_{\min}(z)=-P^{odd}_{\max}(-z).$
\end{theorem}

\begin{proof}
Note, that the function $P^{odd}_{\max}(z)$ is rational with singular points $z_1,z_2,z_3,z_4,$ where $z_3, z_4$ are roots of the equation $1-2z\cos(2\pi/(N+3))+z^2=0,$ and $z_1=-1, z_2=1.$  Computation of the limits at these points indicates that they are removable singularities (see Lemma \ref{le:A11} below). Thus, after the removal of the singularities, the function $P^{odd}_{\max}(z)$ becomes a polynomial of degree $N.$ 

The formula \eqref{maxtrig} can be derived {from the expression for $P^{odd}_{\max}(e^{it})$ by taking the imaginary part}. It implies that ${\rm Im}\left(P^{odd}_{\max}(e^{it})\right)\ge 0$ for $t\in(0,\pi)$. Therefore, $P^{odd}_{\max}(z)$ is a typical real polynomial. 

Applying Taylor's formula we get $P^{odd}_{\max}(z)=z+2\cos(2\pi/(N+3))z^2+ o(z^2).$ Since the coefficient in front of $z^2$ is taking the maximal value, by the uniqueness of the extremal polynomial, we conclude that $P_{max}(z)$ is indeed the required extremizer. 
\end{proof}

\begin{corollary} Straightforward computations lead to formulas
    $$
    P_{max}^{odd}(1)=\frac{\cos(\frac{2\pi}{N + 3})+2}{2-2\cos(\frac{2\pi}{N + 3})},\quad\mbox{and}\quad P_{max}^{odd}(-1)= \frac{\cos(\frac{2\pi}{N + 3})-2}{2+2\cos(\frac{2\pi}{N + 3})}
    $$
\end{corollary}

We note the following:

i) The extremal polynomial is written as a sum of two rational functions such that the first one has only real values on the central unit circle. The imaginary part of the second function on the upper semicircle determines a bounded non-negative trigonometric kernel. \\

ii) Further, note that 
$$
{\rm Im}\left\{\frac{z}{1-2z\sin^2\frac{2\pi}{N+3}+z^2}\right\}=0\quad\mbox{for}\quad z=e^{it}.
$$ 

iii)  Also, observe that the condition Im$\{P^{odd}_{\max}(e^{it})\}=0$ is equivalent to $\sin\frac{(N+3)t}2=0,$ which holds for $t=2\pi/(N+3).$ Since $\cos\frac{2\pi}{N+3}$ is a root of $U_{\frac{N+2}2}(x)$, by Theorem \ref{tm:1}, it is a root of $\Delta_N(x),$ where $\Delta_N(x)$ is the determinant of the matrix ${\bf\Phi}_N(x)$ (see Section \ref{sec:3.2}). Hence, condition Im$\{P^{odd}_{\max}(e^{it})\}=0$ implies $\Delta_N(\cos t)=0.$\\

Now, we apply formulas ~\eqref{eq:a_2j+1}, \eqref{eq:a_2j} 
 for $N=3,$ i.e. compute $P^{odd}_{\max}(z):=z+a_2z^2+a_3z^3$. For the coefficients $a_2$ and $a_3$  we choose $j=1.$ Then, 
$$
a_{2}=\frac23U_0(\mu_3)\left[2 U_1(\mu_3)-\mu_3U_0(\mu_3) \right].
$$
Because $\mu_3=\cos(\pi/3)=1/2,\; U_0(x)=1,$ and $U_1(x)=2x$ we have $a_2=1.$

Similarly,
$$
a_3=\frac13\Bigg[\mu_3U_0(\mu_3)U_1(\mu_3)+2(U_1(\mu_3))^2-(U_0(\mu_3))^2\Bigg]=\frac12.
$$
Hence, $\displaystyle P^{odd}_{\max}(z)=z+z^2+\frac12z^3.$ This result can also directly be achieved by 
Theorem \ref{tm:im-odd} for $N=3$:
\begin{align}\notag
&P^{odd}_{\max}(z)=\frac z{z^2+1-z} + \frac{z^3(-z^6+1)}{2(1-z)(1+z)(z^2+1-z)^2}=z+z^2+\frac12z^3,\\
&P^{odd}_{\min}(z)=-P^{odd}_{\max}(-z)=z-z^2+\frac12z^3.\notag
\end{align}

\section{The case of even $N$}
\subsection{\bf Eigenvectors of the matrix pencil $\bf \{A-\lambda B,\lambda\in\mathbb{C}\}$}
{
\begin{theorem}\label{tm:z1}
Let $\eta$ be a root of the equation $\Delta_N=0$. The solution of the system of linear equations
\[
(2\eta\mathbf B-\mathbf A)\mathbf Z=\mathbf 0
\]
is the one-parameter family $c{\mathbf Z}^{(1)}(\eta)$, where $c\in\mathbb{R}$, ${\mathbf Z}^{(1)}(x)=(z_1^{(1)}(x),\ldots,z_N^{(1)}(x))^T,$ and for $1\le k\le N/2$  
\begin{align*}
&\begin{aligned}
z_{2k-1}^{(1)}(x)=U_{k-1}(x)U_{k-1}(x) - R_N(x)\Big(U_{2k-1}(x)+2k\Big),\\
\end{aligned} \\
&\begin{aligned}
z_{2k}^{(1)}(x)= U_{k-1}(x)U_k(x) - R_N(x)\Big(U_{2k}(x)-2k-1\Big),\\
\end{aligned} \\
& R_N(x)=\frac{(N+2)(N+4)}{4(N+3)}\frac1{q_N^2(x)-1},\qquad q_N(x)=
\frac{N+2}{2U_{N/2}(x)}.
\end{align*}
\end{theorem}
\begin{proof} The proof follows from Lemma \ref{le:A6} and Lemma \ref{le:A7a}.\end{proof}
}

{
In the following, we need the maximal root of the equation $\Delta_N=0$ (see \eqref{eq:deltan} for $\Delta_N$). Denote it by $\eta_N$, and the 
corresponding eigenvector by $\mathbf{Z}^{(1)}(\eta_N).$ If $\nu_N$ is the minimal positive root of the equation 
$U'_{N+2}\left(x\right)=0$ then, by \eqref{eq:U'U'}, $\eta_N=1-2\nu_N^2.$

\subsection{\bf Computing the coefficients of the extremizers}
Let $N$ be an even number. Then, by formulas~\eqref{eq:gamma-0}, for $1\le j\le N/2$ we have
\begin{align}\label{eq:gammas}
\begin{aligned}
\gamma_{2j-1}(x)&=\sum_{k=1}^{N-2j+2} z_k^{(1)}(x) z_{k+2j-2}^{(1)}(x), \qquad
\gamma_{2j}(x)&=\sum_{k=1}^{N-2j+1} z_k^{(1)}(x) z_{k+2j-1}^{(1)}(x);
\end{aligned}
\end{align}
\begin{align}\label{eq:as}
\begin{aligned}
a_{2j-1}=\frac{\gamma_{2j-1}(\eta_N)-\gamma_{2j+1}(\eta_N)}{\gamma_1(\eta_N)-\gamma_3(\eta_N)}, \qquad
a_{2j}=\frac{\gamma_{2j}(\eta_N)-\gamma_{2j+2}(\eta_N)}{\gamma_1(\eta_N)-\gamma_3(\eta_N)},
\end{aligned}
\end{align}
$\gamma_{N+1}=\gamma_{N+2}=0$.
 In the formulas above, we separated the odd and even coefficients to stress that, regardless of the parity of $N$, the odd and even coefficients are computed differently (see Theorems \ref{tm:3} and \ref{tm:z1}.)

\subsubsection{\bf Compact form for extremizers}\label{sec:compform}
{
In the even case, it would be natural to use formulas \eqref{eq:gammas} and \eqref{eq:as}. However, due to the far more complicated formulas obtained from \eqref{eq:gammas} and \eqref{eq:as}, we did not succeed in modifying the approach for odd $N$ to obtain a compact representation. 
Fortunately, the form and the properties of the odd case extremizers allowed us to make an educated guess, which turned out to be correct by verification.
 
\begin{theorem}\label{tm:evencomp} Let $N$ be even. Then, the extremal polynomial allows the following representation
$$P^{even}_{max}(z)=P^e_1(z)+P^e_2(z).$$
Here 
\begin{equation}\label{aseq:3-1}
P^e_1(z)=\frac{z+z^5+\gamma_1(z^2+z^4)+\gamma_2z^3}{(1+z)^2(z^2+1-2\eta_Nz)^2},
\end{equation}
$$
\eta_N=1-2\nu_N^2,\qquad \gamma_1=2(1-\eta_N),\qquad \gamma_2=\frac2{N+3}(-2\eta_N^2-2(N+3)\eta_N+N+5);
$$
\begin{equation}\label{aseq:2-1}
P^e_2(z)=Q_N\frac{2^3z^4}{(1-z)(1+z)^3(z^2+1-2\eta_Nz)^2}
\left(\left(\frac{N+4}2\right)^2(1-z^{N+2})+\right.
\end{equation}
$$
\notag \left.\left(\frac{N+2}2\right)^2(1-z^{N+4})  
+\frac{(N+2)(N+4)}2(z-z^{N+3})\right),
$$
where $\displaystyle Q_N=\frac{2(1-\eta_N^2)}{(N+2)(N+3)(N+4)}$.

For the problem $\min\{a_2\}$, the extremizer is obtained from the maximum problem by alternating signs for even powers.

The resulting non-negative trigonometric polynomial has the form
\begin{equation}\label{aseq:1}
    {\rm Im}\left(P^{even}_{\max}(e^{it})\right)=Q_N\frac1{1+\cos t}
\frac1{\sin t}\frac{\left(\frac{N+4}2\sin\frac{N+2}2t + \frac{N+2}2\sin\frac{N+4}2t\right)^2}{(\cos t-\eta_N)^2}.
\end{equation} 
\end{theorem}

\begin{proof} One can verify that the function $P^{even}_{max}(z)$ does not have poles by standard methods (see Lemma \ref{le:A13} below). If we define this function in the removable singularity points by continuity then it becomes a polynomial of degree $N.$ Formula \eqref{aseq:3-1} implies \eqref{aseq:1}, so that ${\rm Im}\{P^{even}_{\max}(e^{it})\}\ge 0$ for $t\in[0,\pi].$ Thus, the polynomial $P^{even}_{\max}(z)$ is typically real. The Taylor expansion for this polynomial yields $P^{even}_{max}(z)=z+2\eta_N z^2 +o(z^2).$ The coefficient in front of $z^2$ is the maximal possible. 
\end{proof}

\begin{corollary} Straightforward computations lead to formulas
    $$
    P_{max}^{even}(1)=\frac{\mathit{\eta_N}+2}{2-2 \mathit{\eta_N}},\qquad\mbox{and}\qquad
P_{max}^{even}(-1)=\frac{\mathit{\eta_N}-4}{6 \mathit{\eta_N}+6}.
    $$
\end{corollary}

Let us provide some heuristics that lead to \eqref{aseq:3-1} and \eqref{aseq:2-1}.  We conjectured that the structure of the odd representations and the even ones are similar. Because of \eqref{eq:a2-estimate}, it is suggestive to simply replace $\mu_n$ by $\eta_N.$
Note that $\eta_N$ is a maximal root of the equations $\displaystyle U^\prime_{N/2+1}(x)- U^\prime_{N/2}(x)=0,$ and  
$\frac{N+4}2U_{N/2}(x)+ \frac{N+2}2 U_{N/2+1}(x)=0.$ Further, in the odd case the factor $\sin\frac{N+3}2t$ guarantees the property: if Im$\{P^{odd}_{\max}(e^{it})\}=0$ then $\Delta_N(\cos t)=0.$ In the search for a substitute for $\sin\frac{N+3}2t$ providing the corresponding feature, we arrived at the scaled version of the left-hand side of \eqref{eq:RSt} with the normalization factor $Q_N.$ Finally, the choice of $P_1^e(z)$ enabled us to get rid of the poles in  $P_1^e(z)+P_2^e(z).$

Let us mention that the representation of typically real polynomials in the form of a rational function is not new, e.g. it can be found in \cite{MB89}. However, the authors were not able to derive the representation in Theorem \ref{tm:evencomp} from the results in \cite{MB89}.\\ 
}

Finally, let us illustrate formulas \eqref{aseq:3-1} and \eqref{aseq:2-1}. For the case $N=2$, the estimate \eqref{eq:a2-estimate} is $|a_2|\le2(1-2\nu_2^2),$ where $\nu_2$ is the smallest positive root of the equation $U'_4(x)=(16x^4-12x^2)'=0,$ from where $\nu_2^2=3/8.$ Hence $|a_2|\le 1/2$ is a sharp estimate. Further, the equation
$$U^\prime_2(x)- U^\prime_1(x)=(4x^2-1)^\prime-(2x)^\prime=0$$ 
implies the root $\eta_2=\frac14.$ Then, $Q_2=\frac1{64},$ $\gamma_1=\frac32,$ $\gamma_2=\frac74.$
Therefore,
\begin{align}\notag
&P^e_1(z)=\frac{z+z^5+\frac32(z^2+z^4)+\frac74z^3}{(1+z)^2\left(z^2+1-\frac12z\right)^2},\quad P^e_2(z)=\frac{z^4(-4z^6-12z^5-9z^4+12z+13)}{2^3(1-z)(1+z)^3\left(z^2+1-\frac12z\right)^2},\notag\\
&P^{even}_{\max}(z)=P^e_1(z)+P^e_2(z)=z+\frac12z^2,\quad P^{even}_{\min}(z)=z-\frac12z^2.\notag
\end{align}

\section{Summary and Remarks} We can summarize the present results in the following 
\begin{theorem}\label{tm:main}
For typically real polynomials $P(z)=z+\sum_{j=2}^N a_j z^j$ on the unit disk $\mathbb{D},$ the following
exact estimates are valid:
\begin{equation}\notag
|a_2|\le\left\{
\begin{array}{ll}
2\mu_N,\;N\text{ is odd}, \\
2\eta_N=2(1-2\nu_N^2),\;N\text{ is even};
\end{array}
\right.
\end{equation}
where $\mu_N=\cos\frac{2\pi}{N+3}$, $\nu_N$ is the smallest positive root of the equation $U'_{N+2}(x)=0,$ and $\eta_N$ is the maximal root of $\displaystyle U^\prime_{\frac N2+1}(x)- U^\prime_{\frac N2}(x)=0.$ 

In the case of odd $N$, the coefficients of the extremal polynomial for the upper bound are defined by formulas~\eqref{eq:a_2j+1} and \eqref{eq:a_2j}, and for even $N$ by formulas~\eqref{eq:gammas} and \eqref{eq:as}. Concerning the lower bounds, 
the coefficients with even indices are taken with a minus sign.  

The compact form of the extremal polynomials is given by formulas \eqref{ext1} and \eqref{ext} in the odd case, and by \eqref{aseq:3-1} and \eqref{aseq:2-1} in the even one. 
\end{theorem}

\subsection{An estimate for $a_3$}  
Of course, there is also the question of whether the above approach can be used to attack the case {\rm max}$\{a_j\},\; 3\le j\le N-2.$}} Let us look back.
Since we are dealing with optimization in finite dimension, it would be natural to expect that the problem can be reduced to some matrices. A departure point is the matrix representation of the Chebyshev polynomials 
$$
 U_N(x)=\det\begin{pmatrix}
2x & 1 & 0 & 0 & \ldots &0 \\
1 & 2x & 1 & 0 & \ldots &0\\
0 & 1 & 2x & 1 & \ldots &0\\
\ldots & \ldots & \ldots & \ldots & \ldots\\
0 & \ldots & \ldots&0  & 1& 2x \\
\end{pmatrix},
$$
which is a 3-band matrix, or the product representation
 $$
 U_N(x)=2^N \prod_{k=1}^N\left(x-\cos\frac{k\pi}{N+1}\right). $$
 Formally, the above two formulas allow to find eigenvalues of the corresponding matrix. 

This led us to the determinant of the 5-band matrix \eqref{5band}, in which, not only the Chebyshev polynomials are involved, but also their derivatives. In the current article, we are dealing with the 7-band matrix from Lemma \ref{le:A2}. 

Chebyshev polynomials and their derivatives play an essential role. The roots of the polynomials give the extremal values in the odd case, while the roots of the derivatives are in charge in the even case. The roots of the Chebyshev polynomials are much easier to handle in the computations than those of their derivatives. This explains the greater effort concerning the even case in this paper.

Now, a general method for the estimates of the coefficients and a method of finding extremizers becomes clearer. Say, for the $a_3$ coefficient, the corresponding matrix $B$ remains while the matrix $A$ is a modification by inserting a diagonal of zeros above and below the main diagonal in the matrix $A$ from the current paper, i.e.
$$
\begin{pmatrix}
0 & 0 &1/2 & 0 & -1/2 & \ldots \\
0 & 0 & 0 & 1/2 & 0 & \ldots \\
1/2& 0& 0 & 0  & 1/2& \ldots\\
0 & 1/2 & 0 & 0 & 0&\ldots \\
-1/2 & 0 & 1/2 & 0 & 0&\ldots \\
\ldots & \ldots & \ldots & \ldots & \ldots
\end{pmatrix}
$$
Because sharp estimates for $a_3$ are known due to Rogosinski-Szeg\"o \cite{RS50}  and Ruscheweyh \cite{R86}, the natural problem would be to deduce the corresponding extremizers for the $a_3$ case.

\section{Appendix}
\setcounter{theorem}{0}
\renewcommand{\thetheorem}{A.\arabic{theorem}}
\setcounter{equation}{0}
\renewcommand{\theequation}{A.\arabic{equation}}
\begin{lemma}\label{le:A1}
The matrix
\[
\mathbf B=
\begin{pmatrix}
1 & 0 & -1/2 & 0 & \ldots \\
0 & 1 & 0 & -1/2 & \ldots \\
-1/2 & 0 & 1 & 0 & \ldots \\
0 & -1/2 & 0 & 1 & \ldots \\
\ldots & \ldots & \ldots & \ldots & \ldots
\end{pmatrix}
\]
is positive definite; moreover, the successive principal minors $B_k$ $(k=1,\ldots,N)$ of this matrix are defined by the
formula
\begin{align*}
B_k=\left\{
\begin{aligned}
&\frac{(k+2)^2}{2^{k+2}},\; &k\text{ is even},\\
&\frac{(k+1)(k+3)}{2^{k+2}},\; &k\text{ is odd}.
\end{aligned}
\right.
\end{align*}
\end{lemma}
\begin{proof}
Suppose that the matrix $B$ has dimension $N \times N$. Apply the formula~\cite{DSS22,DSS23}
\begin{equation}\label{5band}
  \det
\begin{pmatrix}
1-4x^2 & 2x^2 & -1/2 & 0 & \ldots \\
2x^2 & 1-4x^2 & 2x^2 & -1/2 & \ldots \\
-1/2 & 2x^2 & 1-4x^2 & 2x^2 & \ldots \\
0 & -1/2 & 2x^2 & 1-4x^2 & \ldots \\
\ldots & \ldots & \ldots & \ldots & \ldots
\end{pmatrix}
=\frac{(-1)^N}{2^{N+2}x} U_{N+1}(x) U'_{N+1}(x).  
\end{equation}

Then,
\[
B_k=\frac{(-1)^k}{2^{k+2}} \lim_{x\to0}\frac{U_{k+1}(x) U'_{k+1}(x)}x.
\]

When $k$ is even, we have 
\begin{align*}
&\frac1x U_{k+1}(x)=(-1)^{k/2}(k+2)+\ldots, \quad U'_{k+1}(x)=(-1)^{k/2}(k+2)+\ldots,\\
&\frac1x U_{k+1}(x) U'_{k+1}(x)=(k+2)^2+\ldots,
\end{align*} 
and for odd $k$
\begin{align*}
&U_{k+1}(x)=(-1)^{(k+1)/2}+\ldots, \quad \frac1x U'_{k+1}(x)=-(-1)^{(k+1)/2}(k+1)(k+3)+\ldots, \\
&\frac1x U_{k+1}(x) U'_{k+1}(x)=-(k+1)(k+3)+\ldots.
\end{align*} 
The symbol ``$\ldots$'' denotes the terms containing positive powers of $x$. 

Hence, 
$$
B_k=\frac1{2^{k+2}}\left\{
\begin{array}{ll}
(k+1)(k+3), \;N\text{ is odd}, \\
(k+2)^2,\hspace{1.2cm}N\text{ is even}.
\end{array}
\right.
$$

Overall, all successive principal minors of the matrix $B$ are positive, which yields that this matrix is positive 
definite. The lemma is proved. 
\end{proof}

\begin{lemma}\label{le:A2}
Consider the seven-band matrix
\[
{\bf\Phi}_N (x)=
\begin{pmatrix}
2x & -1/2 & -x & 1/2 & 0 & 0 & \ldots \\
-1/2 & 2x & -1/2 & -x & 1/2 & 0 & \ldots \\
-x & -1/2 & 2x & -1/2 & -x & 1/2 & \ldots \\
1/2 & -x & -1/2 & 2x & -1/2 & -x & \ldots \\
0 & 1/2 & -x & -1/2 & 2x & -1/2 & \ldots \\
0 & 0 & 1/2 & -x & -1/2 & 2x & \ldots \\
\ldots & \ldots & \ldots & \ldots & \ldots & \ldots & \ldots  
\end{pmatrix}
\]
of dimension $N\times N$ $(N>20)$. Let $\Delta_k$, $k=1,\ldots,N$, denote the successive principal minors of this 
matrix $(\Delta_N=\det{\bf\Phi}_N)$. Then, there holds the relation
\begin{equation}\label{eq:le-determ-relation}
\Delta_N-\sum_{j=1}^{10} 2^{-j} b_j\Delta_{N-2j}=0,
\end{equation}
where

\centering
\begin{minipage}{.4\textwidth}
\begin{align*}
b_1&=8x^2-3, \\
b_2&=-24x^4+16x^2-13/4, \\
b_3&=32x^6-24x^4+8x^2-1, \\
b_4&=-16x^8+6x^4-4x^2+7/8, \\
b_5&=16x^8-16x^6+12x^4-5x^2+7/8,
\end{align*}
\end{minipage}
\begin{minipage}{.3\textwidth}
\begin{align*}
b_6&=2^{-2}b_4, \\
b_7&=2^{-4}b_3, \\
b_8&=2^{-6}b_2, \\
b_9&=2^{-8}b_1, \\
b_{10}&=-2^{-10}.
\end{align*}
\end{minipage}
\end{lemma}
\begin{proof}
Formula~\eqref{eq:le-determ-relation} is deduced by repeated application of the Laplace expansion rule for determinants. 

In the expressions for $b_k$, $k=1,\ldots,10$, make the change of variables $x=\frac12(z+z^{-1})$ and denote
$\hat{b}_k=b_k\Big|_{x=\frac12(z+z^{-1})}$, $k=1,\ldots,10$. Then,
\begin{align*}
\hat{b}_1&=1+2(z^2+z^{-2}), \\
\hat{b}_2&=-\frac{17}4-2(z^2+z^{-2})-\frac32(z^4+z^{-4}), \\
\hat{b}_3&=4+\frac72(z^2+z^{-2})+\frac32(z^4+z^{-4})+\frac12(z^6+z^{-6}), \\
\hat{b}_4&=-\frac{13}4-3(z^2+z^{-2})-\frac{11}8(z^4+z^{-4})-\frac12(z^6+z^{-6})-\frac1{16}(z^8+z^{-8}), \\
\hat{b}_5&=\frac94+\frac32(z^2+z^{-2})+(z^4+z^{-4})+\frac14(z^6+z^{-6})+\frac1{16}(z^8+z^{-8}), \\
\hat{b}_6&=2^{-2}\hat{b}_4,\quad \hat{b}_7=2^{-4}\hat{b}_3,\quad \hat{b}_8=2^{-6}\hat{b}_2,\quad 
\hat{b}_9=2^{-8}\hat{b}_1,\quad \hat{b}_{10}=-2^{-10}.
\end{align*}\end{proof}

\begin{lemma}\label{le:A3}
Consider the equation
\begin{equation}\label{eq:psi}
\Psi_n-\sum_{j=1}^{10} 2^j\hat{b}_j\Psi_{n-j}=0.
\end{equation}
The functions $\Psi_n=\alpha+\beta n$, $\Psi_n=(\gamma+n)z^{2n}$, $\Psi_n=(\gamma+n)z^{-2n}$ are particular
solutions of equation~\eqref{eq:psi} ($\alpha$, $\beta$, $\gamma$ are arbitrary constants).
\end{lemma}
\begin{proof}
Substitute the function $\Psi_n=\alpha+\beta n$ into \eqref{eq:psi} and multiply both sides of the equality by $z^8$. On the
left-hand side, we get a polynomial in $z$ of degree $16$. Performing identical transformations, we see that all the 
coefficients of this polynomial are equal to zero. Proceed analogously with $\Psi_n=(\gamma+n)z^{2n}$: multiply by
$z^{20-2n}$, obtain a polynomial in $z$ of degree $20$, and make sure that all the polynomial coefficients equal zero. 
It can be similarly shown for the function $\Psi_n=(\gamma+n)z^{-2n}$. The lemma is proved.    
\end{proof}

\begin{lemma}\label{le:A4}
The following identities hold: 
\begin{align*}
\hspace{-3.1cm}\textnormal{a)}\hspace{3.1cm} U'_k (x)&=\frac1{2(1-x^2)}\Big((k+2) U_{k-1}(x)-kU_{k+1}(x)\Big) \\
&=\frac1{1-x^2}\Big((k+1)U_{k-1}(x)-kxU_k(x)\Big),\\
\hspace{-1.1cm}\textnormal{b)}\hspace{1.1cm} U'_{k+1}(x)-U'_k(x)&=\frac1{1+x}\left((k+2)U_k(x)+(k+1)U_{k+1}(x) \right),\\
U'_{k+1}(x)+U'_k(x)&=\frac1{1-x}\left((k+2)U_k(x)-(k+1)U_{k+1}(x) \right),\\
(U'_{k+1}(x))^2-(U'_k(x))^2&=
\frac1{1-x^2}\Big((k+2)^2(U_k(x))^2-(k+1)^2(U_{k+1}(x))^2\Big), \\
\hspace{-1cm}\textnormal{c) }\hspace{1cm} U'_{k+1}(x)+U'_k(x)&=\frac{\sqrt2}{4\sqrt{1+x}} U'_{2k+2}\Bigg(\sqrt{\frac{1+x}2}\Bigg), \\
\hspace{-1cm}\textnormal{d) }\hspace{1cm} U'_{k+1}(x)-U'_k(x)&=\frac{(-1)^k\sqrt2}{4\sqrt{1-x}} U'_{2k+2}\Bigg(\sqrt{\frac{1-x}2}\Bigg), \\
\textnormal{e) } (U'_{k+1}(x))^2-(U'_k(x))^2&=\frac{(-1)^k}{8\sqrt{1-x^2}} 
U'_{2k+2}\Bigg(\sqrt{\frac{1+x}2}\Bigg) U'_{2k+2}\Bigg(\sqrt{\frac{1-x}2}\Bigg).
\end{align*}
\end{lemma}

\begin{proof}
a) \cite[Lemma~2]{DSS22}; b) {Let us prove the first formula (the second one can be proved in the same way while the third one follows from the preceding two).
\begin{align*}
&\begin{aligned}
&2(1-x^2)\left(U'_{k+1}(x)-U'_k(x)\right) \\
&=(k+3)U_k(x)-(k+1)U_{k+2}(x)-(k+2)U_{k-1}(x)+kU_{k+1}(x)\\
&=(n+3)U_k(x)+kU_{k+1}(x)-(k+2)\left(2xU_k(x)-U_{k+1}(x)\right) \\
&\hspace{6cm} -(k+1)\left(2xU_{k+1}(x)-U_k(x)\right)\\
&=2(1-x)\left((k+2)U_k(x)+(k+1)U_{k+1}(x)\right).
\end{aligned} \\
\end{align*}
}

c) Let $y=\cos{t}$, $x=\cos{2t}$ (i.e., $x=2y^2-1$). Then,
\[
\frac{U_{2k+1}(y)}{2y}=\frac{\sin{(2k+2)t}}{2\cos{t}\sin{t}}=\frac{\sin(k+1)2t}{\sin2t}=U_k(x)
\]
or
\[
\hspace{-1.2cm}U_k(x)=\frac{\sqrt2}{2\sqrt{1+x}} U_{2k+1}\Bigg(\sqrt{\frac{1+x}2}\Bigg).
\]
Applying formula~a), we obtain 
\begin{align}\label{eq:lemmaA4-1}
U'_{k+1}(x)+U'_k(x)=\frac1{2(1-x^2)}\Big[(k+3)&U_k(x)-(k+1)U_{k+2}(x) \nonumber \\
&+(k+2)U_{k-1}(x)-kU_{k+1}(x)\Big].
\end{align}
On the other hand, $U'_{2k+2}(y)=\frac{2y}{2(1-y^2)}\left[2(k+2)U_k(x)-2(k+1)U_{k+1}(x)\right]$, hence
\begin{equation}\label{eq:lemmaA4-2}
\frac1{4y}U'_{2k+2}(y)=\frac1{1-x}\Big[(k+2)U_k(x)-(k+1)U_{k+1}(x)\Big].
\end{equation}
It remains to verify that the right-hand sides of formulas~\eqref{eq:lemmaA4-1} and \eqref{eq:lemmaA4-2} are 
identically equal, which is done by direct calculations; 

d) This can be obtained similarly to c) using the relation 
$$U_k(x)=\frac{(-1)^k\sqrt2}{2\sqrt{1-x}}U_{2k+1}\left(\sqrt{\frac{1-x}2}\right);$$

e) This is a corollary of formulas~c) and d). The lemma is proved. 
\end{proof} 

{
\begin{lemma}\label{le:A9}
Let $n>2$ be an even number, $x_{\min}$ be a minimal positive root of the equation $U^\prime_n(x)=0,$ and $x_{\max}$ a maximal root. Then
$$
i)\qquad x_{\min}<\sin\frac\pi n,
$$
$$
ii)\qquad x_{\max}<\cos\frac\pi n.
$$
\end{lemma}

\begin{proof} $i)$ Without loss of generality, we can assume that $n/2$ is an odd number. Then 
$U^\prime_n(0)>0.$ To show that $U^\prime_n(\sin \pi/n)<0$, use formula a) from Lemma \ref{le:A4}. The sign of the quantity  $U^\prime_n(\sin \pi/n)$ coincides with the sign of the function 
$$\omega(t)=(n+2)\sin nt-n\sin(n+2)t$$ 
if $\cos t=\sin \pi/n,$ i.e., if $t=\pi(n-2)/(2n).$ Therefore, $\sin\frac{\pi n(n-2)}{2n}=0$ since $(n-2)/2$ is an integer.

Further,
$$
\sin\frac{\pi n(n+2)(n-2)}{2n}=\sin\left(\frac{\pi(n-2)}2+\frac{\pi(n-2)}n\right)=\sin\frac{2\pi}n,
$$
because $(n-2)/2$ is an even number.

Therefore, $ \omega\left(\frac{\pi(n-2)}{2n}\right)=-n\sin\frac{2\pi}n<0.$ Thus, the function $\omega(t)$ has a zero in $ \left(0,\frac{\pi(n-2)}{2n}\right).$ This implies that the function $U^\prime_n(x)$ has a zero in $\left(0,\sin\frac\pi n\right).$ 

To prove $ii)$ we show that $U^\prime_n(x)>0$ on $\left(\cos\frac\pi n,1\right).$ It is enough to establish that $\omega(t)>0$ for $t\in(0,\frac\pi n).$ To this end, note that $\omega^\prime(t)=2n(n+2)\sin t\sin(n+1)t$ is positive on $\left(0,\frac\pi{n+1}\right):=I_1,$ and negative on $\left(\frac\pi{n+1},\frac{2\pi}{n+1}\right)=:I_2.$ Hence, $\omega(t)$ is increasing on $I_1$ and decreasing on $I_2.$ Since $\frac\pi{n+1}<\frac\pi n<\frac{2\pi}{n+1}$ and $\omega(\frac\pi n)=n\sin\frac{2\pi}n>0$ we get $\omega(t)>0$ for $t\in(0,\frac\pi n)$.
\end{proof}
}

\begin{lemma}\label{le:A5}
Let $N$ be an odd number, and $\mathbf Z^{(0)}(x)=(z_1^{(0)}(x),\ldots,z_N^{(0)}(x))^T$, where $z_N^{(0)}=1$ and
$$
z_{2k-1}^{(0)}(x)=U_{k-1}(x)U_{k-1}(x),\quad z_{2k}^{(0)}(x)=U_{k-1}(x)U_k(x),\quad k=1,\ldots,\mbox{$\frac{N-1}2$}.
$$
Then, for $N\ge5$,
\begin{align*}
{\bf\Phi}_N&(x)\cdot\mathbf Z^{(0)}(x) 
=U_{\frac{N+1}2}(x)\cdot
\left(0,\ldots,0,\frac12 U_{\frac{N-3}2}(x),-\frac12 U_{\frac{N-5}2}(x),-U_{\frac{N-3}2}(x),U_{\frac{N-5}2}(x)\right)^T.
\end{align*}
\end{lemma}
\begin{proof}
{ Write the product $\Phi_N\cdot \mathbf Z^{(0)}(x)$ coordinatewise. Using \eqref{uj}, by direct calculations, we find that the first $N-4$ coordinates are 
identically equal to zero.} Next, it is convenient to make the change $n=(N-1)/2$. Then the last four coordinates of the product reduce to the form 
$$\frac{z^{n+2}-z^{-n-2}}{z-z^{-1}}\left(-\frac12\right)\frac{z^n-z^{-n}}{z-z^{-1}},\qquad
\frac{z^{n+2}-z^{-n-2}}{z-z^{-1}}\frac12\frac{z^{n-1}-z^{-n+1}}{z-z^{-1}},$$ 
$$\frac{z^{n+2}-z^{-n-2}}{z-z^{-1}}\frac{z^n-z^{-n}}{z-z^{-1}}, \qquad
-\frac{z^{n+2}-z^{-n-2}}{z-z^{-1}}\frac{z^{n-1}-z^{-n+1}}{z-z^{-1}},$$ 
which proves the lemma.    
\end{proof}

\begin{remark}
For $N=3$, we have $\mathbf Z^{(0)}(x)=(1,2x,1)^T$.
\end{remark}

\begin{lemma}\label{le:A6}
Let $N$ be an even number, and $\mathbf Z^{(1)}(x)=(z_1^{(1)}(x),\ldots,z_N^{(1)}(x))^T$, where, for $k=1,\ldots,N/2$,
{
\begin{align*}
&\begin{aligned}
z_{2k-1}^{(1)}=-\frac12(U_{N+2}(x)-&N-3)U_{k-1}(x)U_{k-1}(x) \\
& +\frac12 U_{\frac{N}2+1}(x) U_{\frac{N}2}(x)\left(U_{2k-1}(x)-\frac{N+4}{N+2}\frac{U_{\frac{N}2}(x)}{U_{\frac{N}2+1}(x)}2k\right),
\end{aligned} \\
&\begin{aligned}
z_{2k}^{(1)}=-\frac12(U_{N+2}(x)-&N-3)U_{k-1}(x)U_{k}(x) \\
&\;\;+\frac12 U_{\frac{N}2+1}(x) U_{\frac{N}2}(x)\left(
U_{2k}(x) -(2k+1)\right).
\end{aligned} \\
\end{align*}
}
Then, 
\[
{\bf\Phi}_N(x)\cdot\mathbf  Z^{(1)}(x)=\left(0,\ldots,0,\frac14
\bigg(\frac{(N+4)^2}{N+2}\big(U_{\frac{N}2}(x)\big)^2-(N+2)\big(U_{\frac{N}2+1}(x)\big)^2\bigg)\right)^T.
\]
\end{lemma}
\begin{proof}  
Write the product $\mathbf \Phi_N(x)\cdot\mathbf  Z^{(1)}(x)$ coordinatewise. Using \eqref{uj}, by direct calculations, we see that the first $N-1$ coordinates are identically equal to zero. Next, it is convenient to make the change $n=N/2$. The last coordinate of the product will become
\begin{align*}
\frac1{2(n+1)(1-z^2)^4}\Big(-(n+1)^2 z^{2n+10}-(n+1)^2 z^{-2n-2}+(n+2)^2 z^{2n+4}& \\
+(n+2)^2z^{-2n+4}-(3n^2+10n+9)z^{2n+6}-(3n^2+10n+9)z^{-2n+2}& \\
+(3n^2+8n+6)z^{2n+8}+(3n^2+8n+6)z^{-2n}-2(2n+3)z^2(1-z^2)^2\Big).&
\end{align*}
We can rewrite this expression as 
\begin{align*}
\frac{(n+2)^2(z^{n+2}-z^{-n})^2}{2(n+1)(1-z^2)^2}&-\frac{(n+1)(z^{n+4}-z^{-n})^2}{2z^2(1-z^2)^2} \\
&=\frac14\left(\frac{2(n+2)^2}{n+1}(U_n(x))^2-2(n+1)(U_{n+1}(x))^2\right). 
\end{align*}
Now we substitute  $n=N/2$ back. The lemma is proved.  
\end{proof}

{
\begin{remark}\label{re:1}
For $N=2$, we have $\mathbf Z^{(1)}(x)=(1,4x)^T$.
\end{remark}

\begin{remark}\label{re:2}
If $\eta$ is a root of the equation $\displaystyle U^\prime_{\frac N2+1}(x)-U^\prime_{\frac N2}(x)=0$ then, by Lemma \ref{le:A4},
$$
\frac{(N+4)^2}{N+2}(U_{\frac N2}(\eta))^2-(N+2)(U_{\frac N2+1}(\eta))^2=0.
$$
\end{remark}

\begin{corollary} Let $\eta$ be a root of the equation $\displaystyle U^\prime_{\frac N2+1}(x)-U^\prime_{\frac N2}(x)=0$ and the vector $\mathbf Z^{(1)}(x)$ be given as in Lemma {\rm \ref{le:A6}}. Then,
$$
\mathbf{\Phi}(\eta)\cdot \mathbf Z^{(1)}(\eta)=0.
$$
\end{corollary}
\medskip

Now, let $N$ be an even positive integer and $\displaystyle \hat{\mathbf Z}^{(1)}(x)=(\hat z_1^{(1)}(x),\ldots,\hat z_N^{(1)}(x))^T,$
\begin{align*}
&\begin{aligned}
\hat z_{2k-1}^{(1)}(x)=U_{k-1}(x)U_{k-1}(x) - R_N(U_{2k-1}(x)+2k),\\
\end{aligned} \\
&\begin{aligned}
\hat z_{2k}^{(1)}(x)= U_{k-1}(x)U_k(x) - R_N(U_{2k}(x)-2k-1),\; k=1,\ldots,\frac{N}2,\;\\
\end{aligned} \\
& R_N(x)=\frac{(N+2)(N+4)}{4(N+3)}\frac1{q_N^2(x)-1},\; \qquad q_N(x)=\frac{N+2}{2U_{N/2}(x)}.
\end{align*}

\begin{lemma}\label{le:A7a} If $\eta$ is a root of the equation $\displaystyle U^\prime_{\frac N2+1}(x)-U^\prime_{\frac N2}(x)=0$ then the vectors ${\mathbf Z}^{(1)}(\eta)$ and $\displaystyle \hat{\mathbf Z}^{(1)}(\eta)$ are linearly independent.
\end{lemma}

\begin{proof} Let us show that 
$$
\mathbf Z^{(1)}(\eta)=-\frac12\left(U_{N+2}(\eta)-N-3\right)\mathbf{\hat Z}^{(1)}(\nu).
$$
The condition $\displaystyle U^\prime_{\frac N2+1}(\eta)-U^\prime_{\frac N2}(\eta)=0$ and 
 Lemma \ref{le:A4} implies that 
$$
\frac{N+4}{N+2}\frac{U_{\frac N2}(\eta)}{U_{\frac N2+1}(\eta)}=-1.$$ 
What is left to show is that
$(U_{N+2}(\nu)-N-3)R_N(\eta)=U_{N/2}(\eta)U_{N/2+1}(\eta).$ 

For convenience, denote $u:=U_{N/2}(\eta),\; v:=U_{N/2+1}(\eta).$ Because $U_{n+2}(x)=\left(U_{\frac N2+1}(x)\right)^2-\left(U_{\frac N2}(x)\right)^2,$ the desired identity is reduced to the following,
$$
(v^2-u^2-N-3)\frac{(N+2)(N+4)}{4(N+3)}\frac{u^2}{\frac{(N+2)^2}4-u^2}=uv,
$$
which can be verified by the substitution $v=-\frac{N+4}{N+2}u.$
\end{proof}
}

\begin{lemma}\label{le:A7}
$\displaystyle \sum_{j=0}^n \left(U_j\left(\cos\frac{\pi}{n+2}\right)\right)^2=\frac{n+2}{2\sin^2\frac{\pi}{n+2}}$.
\end{lemma}
\begin{proof}
This follows from the formula
\[
\sum_{j=0}^n \sin^2(j+1)t=\frac14\left(2n+3-\frac{\sin(2n+3)t}{\sin{t}}\right)
\]
at $t=\pi/(n+2)$.
\end{proof}

\begin{lemma}\label{le:A8}
Let $N$ be an odd number and $\mu_N=\cos\frac{2\pi}{N+3}$. The following identities hold: 
\begin{align*}
&\begin{aligned}
\textnormal{a) } 2\sum_{k=1}^{\frac{N-1}2-j} U_k(\mu_N)U_{k+j}(\mu_N)
=\frac1{1-\mu_N^2}&\Bigg[\left(\frac{N-3}2-j\right)\mu_N U_{j-1}(\mu_N) \\
&\hspace{2cm}-\left(\frac{N-1}2-j\right)U_{j-2}(\mu_N)+2\mu_N^2 U_j(\mu_N)\Bigg],
\end{aligned} \\
&\begin{aligned}
\textnormal{b) } 2\sum_{k=1}^{\frac{N-1}2-j} U_{k-1}(\mu_N)U_{k+j}(\mu_N)
=\frac1{1-\mu_N^2}&\Bigg[\left(\frac{N+3}2-j\right)\mu_N U_{j}(\mu_N) \\
&\hspace{4cm}-\left(\frac{N+1}2-j\right)U_{j-1}(\mu_N)\Bigg].
\end{aligned}
\end{align*}
\end{lemma}
\begin{proof}
a) Note that $U_{\frac{N-3}2-j}(\mu_N)=U_{j+1}(\mu_N)$, $U_{\frac{N+1}2-j}(\mu_N)=U_{j-1}(\mu_N)$. Then, 
using the easily verifiable formula
\begin{align*}
2\sum_{k=1}^{n-j}&\sin(k+1)t\sin(k+j+1)t \\
&=(n-j)\cos{jt}-\frac1{2\sin{t}}\sin(2n+3-j)t+\frac1{2\sin{t}}\sin(j+3)t
\end{align*}
together with the formulas
\begin{gather*}
T_j(x)=xU_{j-1}(x)-U_{j-2}(x), \\
\frac12\left(U_j(x)+U_{j+2}(x)\right)=xU_{j+1}(x)=2x^2U_j(x)-xU_{j-1}(x)
\end{gather*}
(where $T_j(x)$ is a Chebyshev polynomial of the first kind), we obtain formula~a).

b) This is proved similarly, taking into account the formulas 
\begin{align*}
&2\sum_{k=1}^{n-j}\sin{kt}\sin(k+j+1)t=(n-j)\cos(j+1)t-\frac1{2\sin{t}}\sin(2n+2-j)t+\frac1{2\sin{t}}\sin(j+2)t,\\
&U_{N-j}(\mu_N)=-U_{j+1}(\mu_N), T_j(x)=xU_{j-1}(x)-U_{j-2}(x),\quad U_{j+1}(x)=2xU_j(x)-U_{j-1}(x).\\ 
\end{align*}
\end{proof}
\bigskip

{
\begin{lemma}\label{le:A11}
Let $N$ be odd and 
$$
R(z):=z(1-z^2)(1+z^2-2yz) + \frac{4(1-y^2)}{N+3}z^3(1-z^{N+3}),
$$
where $y=\cos\alpha,\; \alpha=2\pi/(N+3).$ Then,
$$
i)\; R(\pm1)=0,\quad ii)\; R(e^{\pm i\alpha})=0,\quad iii)\; R^\prime(e^{\pm i\alpha})=0.
$$
\end{lemma}

\begin{proof} Formulas $i)$ and $ii)$ can be easily verified. Let us show that $R^\prime(e^{i\alpha})=0$ and compute
\begin{align*}
R^\prime(z)&=(1-3z^2)(1+z^2-2yz)+2z(1-z^2)(z-y)\\
 &+\frac{4(1-y^2)}{N+3}3z^2(1-z^{N+3}) - 4(1-y^2)z^3z^{N+2},\\
R^\prime(e^{i\alpha})&=2e^{i\alpha}(1-e^{i2\alpha})(e^{i\alpha}-y)-4(1-y^2)e^{i2\alpha}e^{i(N+3)\alpha},
\end{align*}
from where, taking in mind that $y=\cos\alpha,$ we obtain the desired identity. 
\end{proof}

The next lemma is technical and of an auxiliary nature.

\begin{lemma}\label{le:A12}
If $b\sin  at+ a\sin bt=0,$ then
\begin{align*}
&i) \;\;\;\;b^2(1-\cos 2at)+a^2(1-\cos 2bt)+2ab(\cos(a-b)t - \cos(a+b)t)=0,\\
&ii) \;\;\;b^2\sin 2at + a^2\sin 2bt +2ab\sin(a+b)t=0,\\
&iii) \;\;b\sin 2at +a \sin 2bt +(a+b)\sin(a+b)t-(a-b)\sin(a-b)t=0,\\
&iv)\;\; b\cos 2at + a\cos 2bt+(a+b)\cos(a+b)t -(a+b)(1+\cos(a-b)t)=0.
\end{align*}
\end{lemma}

\begin{proof} The formulas follow from the identities
\begin{align*}
&b^2(1-\cos 2at)+a^2(1-\cos 2bt)+2ab(\cos(b-a)t- \cos(a+b)t)\\
&\hspace{7.2cm} =2(b\sin at+ a\sin bt)^2,\\
& b^2\sin 2at+a^2\sin2bt+2ab\sin(a+b)t\\
&\hspace{4cm} =2(b\sin at+a\sin bt)(b\cos at+a\cos bt),\\
& b\sin2at+a\sin2bt+(a+b)\sin(a+b)t-(a-b)\sin(a-b)t\\
& \hspace{4.5cm} =2(b\sin at+a\sin bt)(\cos at +\cos bt),\\
& b\cos 2at + a\cos 2bt +(a+b)\cos(a+b)t - (a+b)(1+\cos(a-b)t)\\
&\hspace{4.8cm} = -2(b\sin at + a\sin bt)(\sin at +\sin bt).
\end{align*}  
\end{proof}

\begin{lemma}\label{le:A13}
Let $N$ be even and
$$
R(z):=\frac1{z^4}p(z)(1-z^2)+\frac{16(1-y^2)}{(N+2)(N+3)(N+4)}\hat p(z),
$$
where 
\begin{align*}
\hat p(z)&=\left(\frac{N+4}2\right)^2(1-z^{N+2})+\left(\frac{N+2}2\right)^2(1-z^{N+4})+
\frac{(N+2)(N+4)}2(z-z^{N+3}),\\
p(z)&=z+z^5+\gamma_1(z^2+z^4)+\gamma_2z^3,\\
\gamma_1&=2(1-y),\quad \gamma_2=\frac2{N+3}(-2y^2-2(N+3)y+N+5),\; y=\cos\alpha,
\end{align*}
where $\alpha$ satisfies
\begin{equation}\label{alpha}
    (N+4)\sin\frac{N+2}2\alpha+(N+2)\sin\frac{N+4}2\alpha=0.
\end{equation} 
Then,
$$
i)\; R(\pm1)=0,\; R^\prime(-1)=0,\; R^{\prime\prime}(-1)=0; \quad ii)\; R(e^{\pm i\alpha})=0,\; R^\prime(e^{\pm i\alpha})=0.
$$
\end{lemma}

\begin{proof} Formula $i)$ can be directly verified. Let us show that $R(e^{i\alpha})=0$ and compute
$$
z^{-3}p(z)\left.\right|_{z=e^{i\alpha}}=2\cos2\alpha+2\gamma_1\cos\alpha+\gamma_2. 
$$
This implies, for $y=\cos\alpha$,
$$
z^{-3}p(z)\left.\right|_{z=e^{i\alpha}}=\frac4{N+3}(1-y^2).
$$
Then,
$$
z^{-4}p(z)(1-z^2)\left.\right|_{z=e^{i\alpha}}=\frac{-8i}{N+3}(1-y^2)\sin\alpha,\qquad \sin\alpha=\sqrt{1-y^2}.
$$
Observe that
$$
\hat p(e^{i\alpha})=\left(\frac{N+4}2\right)^2(1-\cos(N+2)\alpha)+ \left(\frac{N+2}2\right)^2(1-\cos(N+4)\alpha)
$$
$$
 +\frac{(N+2)(N+4)}2(\cos\alpha -\cos(N+3)\alpha)
$$
$$
- i\left[\left(\frac{N+4}2\right)^2\sin(N+2)\alpha + \left(\frac{N+2}2\right)^2\sin(N+4)\alpha \right.
$$
$$
\left. \frac{(N+2)(N+4)}2(-\sin\alpha +\sin(N+3)\alpha)\right]=i \frac{(N+2)(N+4)}2\sin\alpha.
$$
{The last equality uses \eqref{alpha}} and formulas $i)$ and $ii)$ from Lemma \ref{le:A12} with $a=(N+2)/2$ and $b=(N+4)/2.$ Thus, we obtain 
$$
R(e^{i\alpha})=z^{-4}p(z)(1-z^2)\left.\right|_{z=e^{i\alpha}}+
\frac{16(1-y^2)}{(N+2)(N+3)(N+4)}\hat p(e^{i\alpha})=0.
$$

Now, we show that $R^\prime(e^{i\alpha})=0.$ Compute 
$$
\left(z^{-3}p(z)\right)^\prime=z^{-1} \left.\left(2(z^2-z^{-2})+\gamma_1(z-z^{-1}) \right) \right|_{z=e^{i\alpha}}=4ie^{-i\alpha}(1+y)\sin\alpha;
$$
$$
\left(z^{-4}p(z)(1-z^2)\right)^\prime=\left.\left(-\left(z^{-3}p(z)\right)^\prime(z-z^{-1})-z^{-3}p(z)(1+z^{-2})\right)\right|_{z=e^{i\alpha}}
$$
$$
\hspace{1.25cm} =8e^{-i\alpha}(1+y)(1-y^2)-8e^{-i\alpha}\frac{y(1-y^2)}{N+3};
$$
$$
\left(\hat p(z)\right)^\prime=-z^{-1}\frac{(N+2)(N+4)}2\left(\frac{N+4}2z^{N+2} + \frac{N+2}2z^{N+4}-z+(N+3)z^{N+3}\right).
$$
{Apply \eqref{alpha}} and formulas $iii)$ and $iv)$ from Lemma \ref{le:A12} with $a=(N+2)/2$ and $b=(N+4)/2.$ Then,
$$
\frac{16(1-y^2)}{(N+2)(N+3)(N+4)}\left.\left(\hat p(z)\right)^\prime\right|_{z=e^{i\alpha}}=-e^{-i\alpha}\frac{8(1-y^2)}{N+3}(-\cos\alpha +(N+3)(1+\cos\alpha)).
$$
From there, taking into account that $y=\cos\alpha$, we obtain the desired equality, i.e.
$$
R^\prime(e^{i\alpha})=8e^{-i\alpha}(1+y)(1-y^2)-8e^{-i\alpha}\frac{y(1-y^2)}{N+3}
$$
$$
\hspace{3.2cm}-e^{-i\alpha}\frac{8(1-y^2)}{N+3}(-\cos\alpha+(N+3)(1+\cos\alpha))=0.
$$
\end{proof}
}

\bibliographystyle{abbrv}
\bibliography{arxiv}
\end{document}